\documentclass[journal]{new-aiaa}
\usepackage[utf8]{inputenc}
\usepackage{textcomp}
\usepackage{url}
\usepackage{graphicx}
\usepackage{amsmath}
\usepackage{algorithm}
\usepackage{algpseudocode}
\usepackage[version=4]{mhchem}
\usepackage{siunitx}
\usepackage{comment}
\usepackage{longtable,tabularx}
\usepackage{amsfonts}
\usepackage{bm}
\usepackage{booktabs}
\usepackage{rotating}
\usepackage{multirow}
\usepackage{xcolor}

\title{Global Optimization Framework for Automated Low-Thrust Gravity-Assist Trajectory Design}

\author{Ryo Iijima\footnote{Ph.D. Candidate, Department of Aerospace Engineering; also a Visiting Scholar at the School of Aeronautics and Astronautics, Purdue University. Student Member AIAA. Corresponding Author. Email: \texttt{iijima.ryo.p4@dc.tohoku.ac.jp}}}
\affil{Tohoku University, Sendai, Miyagi 980-8579, Japan}

\author{Kenshiro Oguri\footnote{Assistant Professor, School of Aeronautics and Astronautics. Senior Member AIAA.}}
\affil{Purdue University, West Lafayette, Indiana 47907, USA}

\author{Toshinori Kuwahara\footnote{
Professor, Research Center for Green X-Tech, Green Goals Initiative;
Department of Aerospace Engineering, Graduate School of Engineering;
and Research Center for Space Cross-Tech, Green Goals Initiative.
Member AIAA.}}
\affil{Tohoku University, Sendai, Miyagi 980-8579, Japan}

\begin{document}
\maketitle

\begin{abstract}
Gravity-assist trajectory design requires searching numerous leg combinations, but applying nonlinear programming (NLP) for low-thrust trajectory optimization to all combinations is computationally prohibitive. Therefore, conventional frameworks first perform broad search using lightweight trajectory models such as Lambert-based calculation and simplified deep-space maneuver models, while pruning infeasible candidates. However, to maintain efficiency, conventional approaches handle only low-dimensional formulations with limited constraints such as up to a couple of impulse maneuvers per leg. Consequently, they cannot accommodate the constrained multivariable optimization required for low-thrust design, running the risk of prematurely pruning viable trajectories. To address this, this paper introduces a  convex pruning approach capable of solving constrained multivariable problems directly within broad search while retaining the computational efficiency. We then augment the broad search stage with a robust local low-thrust optimization algorithm based on thrust regularization, which together enable global exploration and optimization in an automated fashion. We apply the proposed framework in a BepiColombo-inspired scenario involving nine gravity assists, successfully demonstrating its broad search capability to discover a far greater number of solutions (36,335) compared to a conventional approach (4,664) while identifying a wider feasible launch window by one year.

\end{abstract}

\section{Introduction}
Interplanetary trajectory design involving multiple gravity assists and low-thrust propulsion must determine not only combinations of multiple trajectory legs, but also nonlinear dynamics, lengthy low-thrust control histories~\cite{chai2019review,Benkhoff2021-bc}.  
Therefore, the resulting design space contains strong nonlinearity and a vast number of possible trajectory patterns \cite{mccarty2018parallel}.  
Since optimizing such strongly nonlinear problems is computationally expensive, designers cannot realistically explore the entire design space densely within practical computation time \cite{conway2010spacecraft,englander2017automated}.  
For this reason, practical trajectory design requires a framework that prunes infeasible solutions before applying local optimization to promising candidates \cite{chai2019review,conway2010spacecraft}.

As a representative conventional framework, the Satellite Tour Design Program (STOUR) conducts a broad search in its initial stage to obtain feasible candidates in the design space \cite{Longuski1990-ki}.  
In this stage, STOUR explores combinations of trajectory legs using lightweight trajectory models, such as Lambert-based models, while pruning infeasible candidates \cite{Alemany2007-rx}.  
However, for the broad search of low-thrust trajectories, two critical issues arise if one were to use models that cannot adequately represent the distribution of continuous thrust. 
First, a candidate that appears infeasible due to the excessive $\Delta v$ under a Lambert-based or impulsive model may become feasible when continuous thrust is distributed along the trajectory leg\cite{Petropoulos2004-yd,McConaghy2004-hq,Mereta2018-wz}.   
Second, if the broad search generates an initial candidate that lies too far from the low-thrust feasible region, the local optimizer may fail to converge or may converge to an undesirable local solution~\cite{Petropoulos2000-zf,Wall2009-tv}.

To address this issue, previous studies introduced fast and low-dimensional low-thrust trajectory generation based on shape-based methods into STOUR as STOUR-LTGA \cite{Petropoulos2000-zf,petropoulos2003simple,Petropoulos2004-yd}.  
Shape-based methods are useful because they provide approximate low-thrust initial solutions before local optimization.  
However, shape-based methods may generate unnecessarily large $\Delta v$ because discrepancies can arise between the prescribed geometric model and the actual equations of motion.  
This property becomes problematic in pruning-based broad search because the framework may reject candidates because of geometric-model error rather than true low-thrust infeasibility.

Another representative approach uses a hybrid optimal control problem (HOCP) framework, as implemented in tools such as the Evolutionary Mission Trajectory Generator (EMTG) \cite{englander2017automated}.  
These frameworks directly search for low-thrust trajectories by combining stochastic search methods, such as monotonic basin hopping, with inner-loop trajectory optimization.  
However, local low-thrust optimization in the inner loop requires very high computational cost.  
Therefore, this type of approach cannot practically perform exhaustive or deterministic broad search over the enormous combinatorial design space (in the order of tens of thousands of solutions), and it inherently risks overlooking promising optimal solutions.

To address this challenge within a deterministic framework, Landau  proposed the Efficient Maneuver Placement method in 2018 \cite{Landau2018-yg}.  
This method expands the design space through a fast pruning process that optimizes the trajectory modification caused by a small intermediate impulsive maneuver while satisfying the boundary conditions at the encounter points in a dynamically consistent manner.  
This result demonstrates the effectiveness of introducing intermediate-maneuver trajectory optimization into the broad search stage.  
Landau et al. have also incorporated this method into the deterministic broad search system \textit{Star} \cite{Landau2022-fk}.
On the other hand, Efficient Maneuver Placement achieves high computational efficiency through a Newton-method-based formulation, and this structure makes the method essentially suited to low-dimensional correction problems.
Therefore, it does not readily extend to problems involving multiple design variables and multiple constraints.  
For low-thrust missions, the correction should not be represented as a single impulsive maneuver, but as a bounded continuous-thrust profile distributed along the leg.
Therefore, broad search requires a pruning model that handles low-thrust dynamics while retaining computational efficiency.

To address this requirement, this study proposes convex pruning (CP) for the broad search stage.  
The proposed method linearizes the low-thrust dynamics and solves a convex-programming problem.  
CP enables deterministic solution of multivariable and multi-constraint problems by approximately formulating them as convex optimization problems \cite{Boyd2004-ed,Acikmese2007-gs,Ozaki2022-zz}.  
In contrast to general nonlinear programming, whose convergence and solution quality may depend on the initial guess and the nonconvexity of the problem, convex programming guarantees that any local optimum of the formulated convex problem is also a global optimum, and standard convex optimization problems can be solved deterministically in polynomial time under appropriate assumptions \cite{Boyd2004-ed,Acikmese2007-gs}.
In the proposed framework, we apply CP to near-feasible Lambert legs that conventional pruning would otherwise remove.  
The CP computes a bounded continuous-thrust correction and evaluates whether the leg can satisfy the prescribed terminal conditions under the linearized low-thrust dynamics.  
In addition, when no feasible solution exists, the CP step terminates infeasible cases with very small computation time.  
This property makes the proposed pruning suitable for broad search, where many infeasible candidates exist in the entire design space.

Furthermore, in the subsequent local-optimization phase, conventional low-thrust formulations can suffer from a singularity in the gradient of the control norm, which deteriorates the conditioning of the NLP problem and destabilizes convergence \cite{Oshima2023-gg}. To address this issue, this study introduces thrust regularization, which analytically removes this singularity, into the Sims--Flanagan transcription, a trajectory model commonly used for local optimization after global search \cite{sims1997preliminary,englander2017automated}. This formulation enables the framework to converge the candidates more robustly and obtain optimized solutions.  
We can also pass the low-thrust control history obtained from the CP to the local optimization as an initial control estimate.  
This connection improves the transition from broad search candidate generation to low-thrust trajectory optimization.

The contributions of this paper are twofold.  
First, we embed the CP into a \textit{Star}-based broad search.  
This pruning discovers some legs that are deemed infeasible via impulsive models but are feasible via low-thrust propulsion.  
As a result, the broad search stage can retain low-thrust-feasible trajectory regions that purely impulsive models or endpoint-correction models cannot capture.  
Furthermore, the CP provides not only low-thrust-feasible trajectory legs but also physically meaningful low-thrust control histories.  
We pass these control histories to the subsequent local-optimization stage as initial control estimates, thereby improving the connection between candidate generation in broad search and low-thrust trajectory optimization.  
Second, we add a robust local-optimization process based on a thrust-regularized Sims--Flanagan transcription (SFT) in the subsequent stage, thereby constructing a unified framework for automated low-thrust multi-flyby trajectory design.  
Through the overall flow, the proposed framework enables broader and more robust exploration of the highly nonlinear mission design space.  
We demonstrate the performance of the proposed framework using a BepiColombo-inspired Mercury transfer case involving nine flybys and low-thrust propulsion \cite{Benkhoff2021-bc}.

The remainder of this paper is organized as follows.  
Section~\ref{sec:background} reviews the background of this study, including the dynamical model and the \textit{Star} overview.  
Section~\ref{sec:method} presents the proposed framework, with emphasis on the CP embedded in the \textit{Star}-based broad search.  
This section also describes the local optimization model, including the thrust-regularized SFT used for local optimization.  
Section~\ref{sec:example} applies the proposed framework to a representative Mercury transfer mission and evaluates its performance. Section~\ref{sec:discussion} discusses the numerical implications, methodological contribution, and future work. Finally, Section ~\ref{sec:conclusion} concludes the paper.


\section{Background} \label{sec:background}

\subsection{Dynamical Model}
This study expresses the general equations of motion for the controlled two-body problem as follows, where $\bm{x}$ denotes the orbital state vector and $t$ denotes the time~\cite{battin1999introduction}:
\begin{equation}
\dot{\bm{x}} = \bm{f} (\bm{x},m,\bm{u},t)=
\bm{f}_\mathrm{kep} (\bm{x}) + \bm{F}(\bm{x})\bm{a}_\mathrm{acc}(m,\bm{u},t) \label{equ: dynamical}
\end{equation} 
Here, $\bm{f}_\mathrm{kep}(\bm{x})$ represents the equations of motion for the Keplerian dynamics, and $\bm{F}(\bm{x})$ denotes the control matrix.
In this study, $\bm{a}_\mathrm{acc}$ represents the acceleration generated by low-thrust propulsion. $m$ is the spacecraft mass, and $\bm{u}$ denotes the dimensionless throttle-direction vector. Unless otherwise specified, $\|\cdot\|$ denotes the Euclidean norm. The thrust-induced acceleration ($\bm{a}_\mathrm{acc}$) and the spacecraft mass decrease due to propellant consumption are given by \cite{Maggi2022-ck}:
\begin{align}
\bm{a}_\mathrm{acc}(m,\bm{u},t) = \frac{D N_\mathrm{EP} T_\mathrm{thrust}(t)  \bm{u}}{m}, \quad \dot{m}(\bm{u},t) = -\frac{D N_\mathrm{EP} T_\mathrm{thrust}(t) \| \bm{u}\|}{g_0 I_\mathrm{sp}(t)}, \quad \| \bm{u}\| \le 1 \label{equ: dynamical control}
\end{align}
where $D$ is the duty cycle, $N_{\mathrm{EP}}$ is the number of operating electric-propulsion units, $T_{\mathrm{thrust}}(t)$ is the available thrust magnitude at time $t$, $g_0$ is the standard acceleration due to gravity, and $I_{\mathrm{sp}}(t)$ is the specific impulse of the thruster.
The convex-pruning formulation in Section~\ref{sec:cp} treats
$T_{\mathrm{thrust}}(t)$ and $I_{\mathrm{sp}}(t)$ as constant parameters,
whereas the local-optimization formulation in
Section~\ref{sec:local optimization} evaluates them as functions of
the power supplied to the propulsion system, as defined in
Section~\ref{sec:interpolation}.
In addition, in this paper, we generally express the equations in Cartesian coordinates. We use Modified Equinoctial Elements (MEE) only in the
CP formulation presented in Section~\ref{sec:cp},
where their definitions and corresponding equations of motion
are introduced. Quantities expressed in Cartesian coordinates and MEE are denoted by the subscripts $\mathcal{C}$ and $\mathcal{M}$, respectively. In Cartesian coordinates, the equations of motion in Eq.~\eqref{equ: dynamical} are defined as
\begin{equation}
\bm{x}_\mathcal{C} =
\begin{bmatrix}
\bm{r} \\
\bm{v}
\end{bmatrix}, \quad
\bm{f}_{\mathrm{kep},\mathcal{C}}(\bm{x}_\mathcal{C}) =
\begin{bmatrix}
\bm{v} \\
-\dfrac{\mu}{|\bm{r}|^3}\bm{r}
\end{bmatrix}, \quad
\bm{F}_\mathcal{C}(\bm{x}_\mathcal{C}) =
\begin{bmatrix}
\bm{0}_{3 \times 3} \\
\bm{I}_{3 \times 3}
\end{bmatrix}
\label{equ:kepler}
\end{equation}

where $\bm{r}$ and $\bm{v}$ denote the spacecraft position and velocity vectors, respectively, and $\mu$ is the gravitational parameter of the central body. For the MEE representation, the control matrix, $\bm{F}_\mathcal{M}(\bm{x}_\mathcal{M})$, is given by the Gauss
variational equations.

\subsection{Gravity-Assist Model} \label{sec:gravity assist}

This study models gravity assist using the zero-sphere-of-influence (ZSOI) patched-conic approximation to reduce the computational cost in the optimization \cite{Vavrina2016-fy}. First, the formulation must ensure that the flyby occurs above a prescribed minimum altitude around the celestial body. Therefore, the periapsis altitude is approximately computed from the direction components and magnitudes of the incoming and outgoing velocity vectors, and the condition that this altitude exceeds the limiting altitude is imposed as follows \cite{battin1999introduction}.
\begin{align}
&\theta=\arccos \left({\frac{{\bm{v}^-_\infty} \cdot {\bm{v}^+_\infty}}{{\|\bm{v}^-_\infty\|} \|\bm{v}^+_\infty\|}} \right) \label{equ:delta}, \\
&r_\mathrm{min} - \frac{\mu }{v_\infty^2}\left[\frac{1}{\sin(\frac{\delta}{2})}-1 \right] \leq 0, \quad
r_\mathrm{min} = r_\mathrm{body} +a_{\mathrm{min}}\label{equ:flyby} \\
& \beta = \max(0,\theta - \delta) 
\label{equ:beta}
\end{align}
where $\theta$ is the angle between the incoming and outgoing hyperbolic-excess velocity vectors, $r_{\min}$ is the minimum allowable periapsis radius,
$r_{\mathrm{body}}$ is the radius of the flyby body,
$a_{\min}$ is the minimum allowable altitude above the body's surface,
and $\mu$ is the gravitational parameter of the flyby body. The superscripts $-$ and $+$ denote quantities
immediately before and after the gravity assist, respectively.
Accordingly, $\bm{v}_{\infty}^{-}$ and $\bm{v}_{\infty}^{+}$ are the
incoming and outgoing hyperbolic excess velocity vectors, respectively.  The angle $\delta$ denotes the turning angle achievable
through the planetary gravity assist, while
$\beta=\max(0,\theta-\delta)$ represents the  turning-angle
violation.
The local-optimization formulation assumes an unpowered ballistic
flyby. Under this assumption, the required turning angle determined by
the incoming and outgoing hyperbolic excess velocity vectors equals the
gravity-assist turning angle, such that $\delta=\theta$. 
The flyby inequality can be rewritten in a numerically more convenient form by eliminating the reciprocal dependence on $\sin(\theta/2)$ and $v_\infty^2$, where
$v_\infty=\min(\|\bm{v}_\infty^{-}\|,\|\bm{v}_\infty^{+}\|)$. Using $\sin^2(\theta/2)=(1-\cos\theta)/2$ and $\cos\theta=(\bm{v}_\infty^{+}\cdot\bm{v}_\infty^{-})/v_\infty^2$, Eq.~\eqref{equ:flyby} is transformed into
\begin{equation}
\left(
\frac{r_\mathrm{min} v_\infty^2 }{\mu}+1
\right)^2
\left(v_\infty^2-\bm{v}_\infty^{+}\cdot\bm{v}_\infty^{-}\right)
-2v_\infty^2 \le 0
\label{equ:flyby-constraint-regularized}
\end{equation}
In the local-optimization formulation described in
Section~\ref{sec:local optimization}, this formulation reduces the magnitude of the constraint Jacobian and improves the numerical scaling of the flyby constraints in the local optimization. Although the exact values depend on the specific trajectory case, the reformulated constraints generally keep the Jacobian entries smaller and suppress excessively large sensitivities, including those associated with the intermediate patch constraints in the SFT.
In addition, the local-optimization formulation assumes unpowered flybys and enforces this assumption by equating the magnitudes of the incoming and outgoing hyperbolic excess velocities:
\begin{align}
v_\infty = \| \bm{v}^-_\infty \| = \| \bm{v}^+_\infty \| \label{equ:velocity manitude}
\end{align}

\subsection{Ephemeris Model, Launch Constraint Model, and Spacecraft Propulsion Model} \label{sec:interpolation}

In the local optimization performed in this study, automatic differentiation provides all gradients of the computational model, including those associated with the ephemeris model. However, direct ephemeris evaluations from SPICE do not directly support the construction of an automatic-differentiation computational graph. To address this issue, a differentiable ephemeris model is constructed from SPICE data using interpolation. Cubic spline interpolation, such as that used in the FIRE system, can provide fast ephemeris output \cite{Arora2010-el}. However, the branching operations and piecewise interpolation inherent in this method can complicate the symbolic computational graph and reduce the efficiency and robustness of automatic differentiation in gradient-based optimization. Therefore, this study employs Chebyshev polynomial interpolation, where all coefficients are stored in a single coefficient matrix, to avoid branching structures in the ephemeris evaluation. The implementation first obtains ephemeris data over the search time interval and computes the Chebyshev coefficients for the sampled data. It then stores these coefficients in the coefficient matrix. During evaluation, the model computes the ephemeris state as a matrix-vector product between the coefficient matrix and the Chebyshev basis vector. The interpolation parameters are tuned so that the position error remains below 1 m and the velocity error remains below 0.1 m/s.

To formulate the optimization problem presented, this study introduces empirical models for launch vehicle performance, SEP engine characteristics, and solar power generation. 
Table~\ref{tab:model_coefficients} presents the values for the polynomial coefficients of the launch-vehicle and SEP-engine models, while Table~\ref{tab:input_formal} summarizes the input parameters used for broad search and local optimization.

\begin{table}[!htbp]
\centering
\caption{Polynomial coefficients for the representative models of launch-vehicle (Falcon Heavy Recovery) and solar electric propulsion (QinetiQ T6) }
\label{tab:model_coefficients}
\begin{tabular}{lrrrrrr}
\hline
\textbf{Model}
& $\eta_1$ & $\eta_2$ & $\eta_3$ & $\eta_4$ & $\eta_5$ & $\eta_6$ \\
\hline

$m_{\mathrm{delivered}}(C_3)$
& $-6.05318{\times}10^{-6}$
& $9.71071{\times}10^{-4}$
& $-6.08631{\times}10^{-2}$
& $2.64628$
& $-177.347$
& $6706.42$ \\

$T_{\mathrm{thrust}}(P)$
& $-8.70819{\times}10^{-17}$
& $1.41849{\times}10^{-12}$
& $-8.68767{\times}10^{-09}$
& $2.45217{\times}10^{-05}$
& $0$
& -- \\

$I_{\mathrm{sp}}(P)$
& $-7.03998{\times}10^{-5}$
& $7.23576{\times}10^{-1}$
& $2.22634{\times}10^{3}$
& $0$
& $0$
& -- \\

\hline
\end{tabular}
\end{table}

The launch constraint uses an inverse correlation curve between $C_3$ and deliverable mass, based on the capability of a real launch vehicle. This curve defines the upper bound on the allowable launch performance. To construct the curve, the model interpolates launch vehicle performance data and applies the specific values for the coefficients ($\eta_{\mathrm{LV},1}, \dots, \eta_{\mathrm{LV},6}$) and  a user-defined launch performance margin ($\sigma_{\mathrm{LV}}$). The source data for this curve come from NASA's Launch Vehicle Performance website\cite{NASA_LVPerf}.
The polynomial approximation represents the interpolated data in the following form\cite{Englander2017-ls}:
\begin{align}
m_{\mathrm{delivered}}(C_3) = 
(1-\sigma_{\mathrm{LV}})
\left(
\eta_{\mathrm{LV},1} C_3^5 
\eta_{\mathrm{LV},2} C_3^4 
+ \eta_{\mathrm{LV},3} C_3^3 
+ \eta_{\mathrm{LV},4} C_3^2 
+ \eta_{\mathrm{LV},5} C_3 
+ \eta_{\mathrm{LV},6}
\right), \quad C_3 = v_\infty^2 \label{equ:mass delivered}
\end{align}

This study assumes a solar electric propulsion (SEP) system for the spacecraft and models both the engine performance as a function of supplied power and the power-supply process from the Sun-dependent available power to the power delivered to the propulsion system. First, the model constructs correlation curves for the thrust ($T$) and specific impulse ($I_{\mathrm{sp}}$) as functions of the supplied power from experimental data for a specific engine. Following previous studies~\cite{englander2017automated,Topputo2021-vi}, this study interpolates these data using the following fourth-order polynomials and introduces them into the optimization as differentiable functions:
\begin{align}
T_\mathrm{thrust}(P) &= \eta_{T,1} P^4 + \eta_{T,2} P^3 + \eta_{T,3} P^2 + \eta_{T,4} P + \eta_{T,5}  \label{equ:thrust}\\
I_{\mathrm{sp}}(P) &= \eta_{I_{\mathrm{sp}},1} P^4 + \eta_{I_{\mathrm{sp}},2} P^3 + \eta_{I_{\mathrm{sp}},3} P^2 + \eta_{I_{\mathrm{sp}},4} P + \eta_{I_{\mathrm{sp}},5} \label{equ:isp}
\end{align}
where $\eta_{\mathrm{T},1},\dots,\eta_{\mathrm{T},5}$ are the polynomial
coefficients for the thrust model, and
$\eta_{I_{\mathrm{sp}},1}, \dots, \eta_{I_{\mathrm{sp}},5}$ are those for the specific-impulse model. The values of these coefficients are summarized in Table~\ref{tab:model_coefficients}. The power model computes the supplied power ($P$) by subtracting the spacecraft bus power consumption from the generated power and then applying a user-defined power margin ($\delta_{\mathrm{power}}$) as follows:
\begin{align}
P(r,t) = (1 - \delta_{\mathrm{power}})
P_{\mathrm{generated}}(r,t)
\end{align}
 For solar power generation, the generated power depends on the distance from the Sun. In general, it is inversely proportional to the square of the heliocentric distance. As a simple model, Topputo et al. introduced the following formulation, where $P_0$ represents the generated power at 1 astronomical unit (AU) \cite{Topputo2021-vi}:
\begin{align}
P_{\mathrm{generated}} (r,t) = \frac{P_0 (t)}{r^2}
\end{align}
  Here, $r$ is expressed in AU. The power model also accounts for time-dependent degradation of the solar arrays. Let $\tau_\mathrm{deg}$ denote the degradation rate, and let $P_{\mathrm{0\mathrm{-}BOL}}$ denote the beginning-of-life (BOL) power at departure. The model computes the reference power as
\begin{align}
P_0 (t) = P_{\mathrm{0\mathrm{-}BOL}} (1 - \tau_\mathrm{deg})^{\frac{t - t_1}{365\times24\times3600}} ,
\end{align}
where $t_1$ denotes the launch epoch in seconds, so that $t-t_1$ represents the elapsed time since launch. This formulation shows that, in solar electric propulsion, the available power scales with the inverse square of the distance from the Sun. Since the thrust $T$ is proportional to $P$, the Sun-spacecraft distance affects the thrust in the same manner.


\subsection{\textit{Star}}

\textit{Star} is a mission design tool developed by Landau et al. at the Jet Propulsion Laboratory (JPL) for the global search of patched-conic trajectories \cite{Landau2022-fk}. Based on user-defined constraints, such as time windows and allowable $\Delta v$, \textit{Star} progressively prunes candidates that do not satisfy the prescribed conditions and rapidly narrows the enormous design space to feasible trajectory candidates. \textit{Star} supports variable-sequence searches, recovery of near-feasible solutions through intermediate-maneuver optimization based on Efficient Maneuver Placement, and searches for resonant trajectories \cite{Landau2018-yg}. These capabilities substantially expand the explored design space and allow the tool to retain candidates that conventional pruning would otherwise remove. Although \textit{Star} cannot directly handle trajectories in extremely weak-gravity environments or trajectories driven by low-thrust propulsion, its comprehensive search can provide high-quality initial guesses for subsequent local optimization, where the model accounts for these higher-fidelity physical effects.

\section{Proposed Optimization Framework} \label{sec:method}
\subsection{System Architecture}
\begin{figure}
    \centering
    \includegraphics[width=1.0\linewidth]{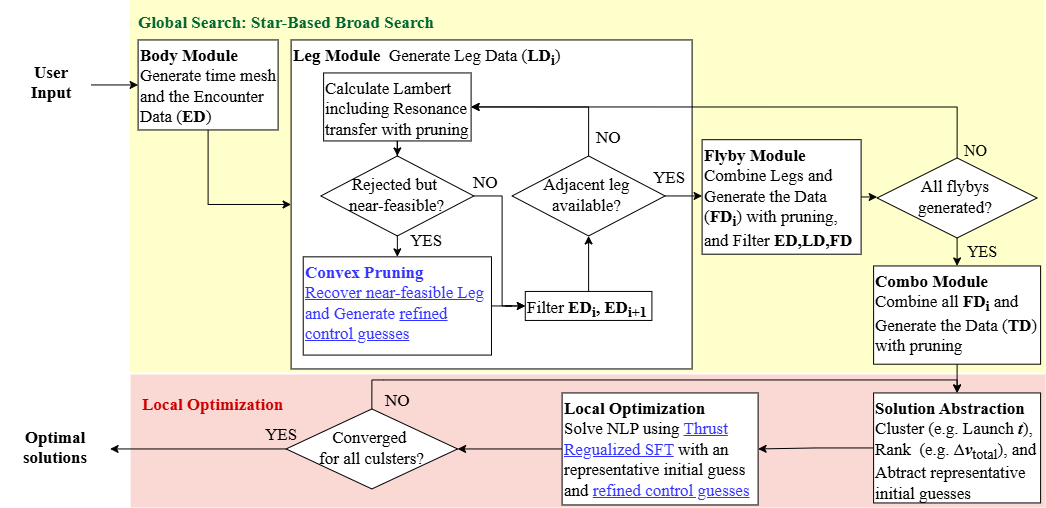}
    \caption{System flow chart of the proposed trajectory design framework integrating CP and thrust-regularized local optimization.}
    \label{fig:SystemFlowChart}
\end{figure}

The proposed framework integrates CP into the \textit{Star}-based broad search, abstracts the resulting candidates by clustering and ranking, and then connects representative candidates to thrust-regularized local optimization, as shown in Fig.~\ref{fig:SystemFlowChart}. The upper part of the figure corresponds to the broad search stage, where trajectory-leg combinations are generated and pruned, while the lower part corresponds to the local-optimization stage, where representative candidates are refined using a higher-fidelity low-thrust model. The key features of the proposed architecture are the CP inserted into the \textit{Leg Module} and the thrust-regularized SFT used in the local-optimization stage; together, these components allow the framework to recover near-feasible Lambert legs before pruning and to robustly refine the selected candidates afterward.

The broad search algorithm is based on the \textit{Star} system \cite{Landau2022-fk}. The user first defines the input parameters, including the sequence of celestial bodies to be visited and the time window for each encounter. Details of the input parameters used in the original \textit{Star} formulation can be found in the previous literature \cite{Landau2022-fk}.
In addition to these parameters, the present framework requires the user to define parameters associated with the CP introduced in this study, such as the maximum thrust magnitude and the near-feasible range in which the pruning is applied.

The broad search consists of four modules: the \textit{Body Module}, \textit{Leg Module}, \textit{Flyby Module}, and \textit{Combo Module}.
Each module generates its corresponding data structure and progressively prunes infeasible candidates according to the user-defined constraints.
In particular, the \textit{Leg Module} generates trajectory legs that connect adjacent encounter points.
In this module, the framework first generates each candidate leg by solving a Lambert problem.
It then applies the proposed CP only when the Lambert-generated candidate is infeasible, is not identified as a resonant-transfer case, and remains within the user-defined recovery range described in Section~\ref{sec:cp}.
By applying the pruning only to candidates that actually require low-thrust correction, the framework reduces the additional computational cost.
Furthermore, this selective application avoids unnecessary thrust usage on trajectory legs that are already feasible or are not appropriate for low-thrust correction, thereby reducing the deviation of the generated initial guess from the low-thrust optimal solution.

After the broad search, the framework clusters the obtained trajectory candidates according to a specified quantity and ranks them within each cluster.
In the numerical example presented in Section~\ref{sec:example}, the framework clusters the candidates by launch epoch, where each cluster corresponds to a one-week launch-time interval.
Within each cluster, it ranks the candidates in ascending order of the total velocity increment, $\Delta v_{\mathrm{total}}$.
This clustering and ranking process corresponds to the solution-abstraction step shown in Fig.~\ref{fig:SystemFlowChart}.

Finally, the framework performs local optimization using the candidate with the smallest $\Delta v_{\mathrm{total}}$ in each cluster as the representative solution.
The local optimization employs the thrust-regularized SFT introduced in Section~\ref{sec:local optimization}.
For gravity assists, the model uses an approximate formulation compatible with the zero-sphere-of-influence (ZSOI) assumption.
The resulting nonlinear programming problem is solved using IPOPT, and automatic differentiation with CasADi provides the derivatives of the computational model \cite{Andersson2018,Wachter2006-io}.
The representative solution selected from the broad search provides the initial trajectory guess.
When the CP provides a low-thrust initial solution, the corresponding thrust history serves as the initial control guess.
Otherwise, the control variables are initialized with a small constant value of $1 \times 10^{-7}$.

After attempting local optimization for the representative solution of each cluster, the framework processes the non-converged clusters again.
For each non-converged cluster, it selects the candidate with the next smallest $\Delta v_{\mathrm{total}}$ as a new representative solution and performs local optimization again.
This process is repeated until a converged solution is obtained for each cluster or until all candidates in the cluster have been tested.
If no candidate in a cluster converges, the framework regards the launch-time region represented by that cluster as containing no feasible solution within the prescribed design space.

\subsection{\textit{Star}-based broad search}

The user defines the mission through several parameters, such as the sequence of celestial bodies and mission constraints. In the \textit{Star} framework, three types of elements describe the trajectory: encounters, legs, and flybys. The departure body has index $1$, and the final destination has index $n_i$.
\begin{equation}
\begin{aligned}
\mathrm{Encounter: } & i = 1,2,\ldots,n_i \\
\mathrm{Leg: } & i = 1,2,\ldots,n_i-1 \\
\mathrm{Flyby: } & i = 2,\ldots,n_i-1
\end{aligned}
\end{equation}
The \textit{Star}-based broad search organizes its computation into four modules: the \textit{Body Module}, \textit{Leg Module}, \textit{Flyby Module}, and \textit{Combo Module}. Each module plays a specific role in the trajectory generation process and stores the intermediate solutions in structured databases with predefined data types.

\begin{algorithm} [t]
\footnotesize
\caption{Broad Search Framework with Convex Pruning}
\label{alg:broad_search}
\begin{algorithmic}[1]
\Require User configurations ${I}_{\mathrm{user}}$ (Time bounds, TOF bounds, $v_{\infty, \max}$, $N_{\mathrm{rev}}$, $\zeta$, $u_\mathrm{max}$, $\Delta v_{\mathrm{leg, max}}$, $\Delta v_{\mathrm{flyby, max}}$)
\Ensure Feasible Trajectory Database ($\mathrm{TD}$)

\State $\mathrm{ED} \gets \textbf{BodyModule}({I}_{\mathrm{user}})$ \Comment{Generate Encounter Data}
\For{$i = 1$ \textbf{to} $n-1$}
    \State $\mathrm{LD}_i \gets \emptyset$
    \For{each valid pair $(p, q)$ in $(\mathrm{ED}_i, \mathrm{ED}_{i+1})$}
        \For{each $N_{\mathrm{rev}} \in {I}_{\mathrm{user}}$}
            \State $[\bm{v}_{\mathrm{D}}, \bm{v}_{\mathrm{A}}] \gets \textbf{LambertSolver}(r_{\mathrm{D}}, r_{\mathrm{A}}, \mathrm{TOF}, N_{\mathrm{rev}})$
            \State Calculate $\bm{v}_{\infty \mathrm{D}}$ and $\bm{v}_{\infty \mathrm{A}}$
            \State $\Delta\theta \gets
            \begin{cases}
            2\pi N_{\mathrm{rev}}+\theta, & \mathrm{short}\\
            2\pi N_{\mathrm{rev}}+(2\pi-\theta), & \mathrm{long}
            \end{cases}$
            
            \If{$\exists N \in \mathbb{Z}^+
            \mathrm{\ s.t.\ } |\Delta\theta-N\pi|<4^\circ$}
                \State $[\bm{v}_{\infty\mathrm{D}}, \bm{v}_{\infty\mathrm{A}}] \gets \textbf{ResonantTransfer}(\bm{r}_{\mathrm{D}}, \bm{v}_{\mathrm{D}}, I_\mathrm{user}, \bm{r}_{\mathrm{A}}, \bm{v}_{\mathrm{A}}, \mathrm{ED}_i, \mathrm{ED}_{i+1})$
            \EndIf

            \If{$(||\bm{v}_{\infty D}|| \le \bm{v}_{\infty,\max,i} )$ and  $(||\bm{v}_{\infty A}|| \le \bm{v}_{\infty,\max,i+1})$} \Comment{Lambert or Resonant Trajectory Accepted}
                \State Add Leg to $\mathrm{LD}_i$
            \ElsIf{(Not Resonant) and ($\max(||\bm{v}_{\infty \mathrm{D}}||, ||\bm{v}_{\infty \mathrm{A}}||) \le \zeta \times \bm{v}_{\infty,\max}$)}
                \State \Comment{\textbf{Proposed Convex Pruning}}
                    \State $[\bm{v}_{\infty \mathrm{D}}^{\mathrm{cp}}, \bm{v}_{\infty \mathrm{A}}^{\mathrm{cp}}, \Delta v_{\mathrm{leg}}] \gets \textbf{ConvexPruning}({I}_{\mathrm{user}},\bm{r}_{\mathrm{D}}, \bm{v}_{\mathrm{D}}, \bm{r}_{\mathrm{A}}, \bm{v}_{\mathrm{A}}, \mathrm{TOF}, \mathrm{ED}_i, \mathrm{ED}_{i+1})$
                    \If{(CP converges) and ($\Delta v_{\mathrm{leg}} \le \Delta v_{\mathrm{leg, max}}$)}
                        \State Add refined Leg($\bm{v}_{\infty \mathrm{D}}^{\mathrm{cp}}, \bm{v}_{\infty \mathrm{A}}^{\mathrm{cp}}$) to $\mathrm{LD}_i$
                    \EndIf
            \EndIf

        \EndFor
    \EndFor
    
    \State Apply $\textbf{EncounterFilter}(\mathrm{ED}, \mathrm{LD}_i)$
    
    \If{$i > 1$}
        \State $\mathrm{FD}_i \gets \textbf{FlybyModule}({I}_{\mathrm{user}},\mathrm{ED}, \mathrm{LD}, \Delta v_{\mathrm{flyby, max}})$
        \State Apply $\textbf{Incoming/OutgoingLegFilters}(\mathrm{ED}, \mathrm{LD}, \mathrm{FD})$
    \EndIf
\EndFor

\State $\mathrm{TD} \gets \textbf{ComboModule}({I}_{\mathrm{user}}, \mathrm{ED}, \mathrm{LD}, \mathrm{FD})$
\State \Return $\mathrm{TD}$
\end{algorithmic}
\end{algorithm}

\subsubsection{Body Module}

Using the time constraints defined by the user, the \textit{Body Module} formulates and solves a linear programming problem to determine the bounds of the encounter times and the time-of-flight (TOF), denoted by $t$ and $\tau$ in seconds, for each leg. First, it computes the encounter-time bounds as
\begin{equation}
\bm{t}_{\min}^* = \arg\min_{\bm{t}} \sum_i t_i, 
\quad
\bm{t}_{\max}^* = \arg\max_{\bm{t}} \sum_i t_i
\end{equation}
subject to
\begin{align}
\bm{t}_{\min} &\leq \bm{t} \leq \bm{t}_{\max},\\
\bm{\tau}_{\min} &\leq \bm{t}A_\mathrm{LP} \leq \bm{\tau}_{\max}.
\end{align}
Here, the matrix $A_\mathrm{LP}$ maps the encounter-time variables $\bm{t} = [t_1, \ldots,t_{n_i}]$ to the corresponding TOF constraints $\bm{\tau}$. 
Second, using the obtained encounter-time bounds, the \textit{Body Module} computes the TOF bounds as
\begin{equation}
\bm{\tau}_{i}^* = \arg\max_{\bm{\tau}} \sum_{k=1}^{n_i} \tau_{i,k}, 
\end{equation}
subject to
\begin{equation}
t_{\min,q} - t_{\max,p}
\leq
\tau_{i,q} - \tau_{i,p}
\leq
t_{\max,q} - t_{\min,p} , \quad \forall i \in \mathbb{Z}_1^{n_i-1}
\end{equation}
The formulation defines the TOF bounds not only for adjacent encounters but also for all encounter pairs. The indices $p$ and $q$ denote the encounter indices satisfying $p < q$. $\mathbb{Z}_a^b$ denotes the set of integers from $a$ to $b$.

Next, the \textit{Body Module} generates encounter data (ED) for each encounter based on the defined time bounds and the number of discretization steps. For each encounter, the module generates a time mesh independently, without using data from other encounters. Therefore, the number of data points can differ among encounters. 
The identifier $I_{\mathrm{E}}$ is assigned to each data entry in the \textit{Body Module}, while
$E$ identifies the celestial body associated with that entry.
In the \textit{Star} framework, the model does not consider the sphere of influence and treats celestial bodies as point locations. Therefore, the departure and arrival positions share the same position vector $\bm{r}$.

\begin{equation}
\begin{aligned}
&\mathrm{ED}_{i,j} = [I_{\mathrm{E},i,j},\, t_{i,j},\, E_{i,j},\,  \bm{r}_{i,j},\, \mu_{i,j},\,  r_{\min,i,j}], \\
&\text{Number of encounters} :\quad i = 1,\ldots,n_i, \\
&\text{Number of data in each encounter} :\quad j = 1,\ldots,n_j
\end{aligned}
\end{equation}

\subsubsection{Leg Module}

In the \textit{Leg Module}, the module generates leg data (LD) for each trajectory segment by solving Lambert problems using adjacent encounter data obtained from the \textit{Body Module}. Each LD entry stores the encounter IDs of the adjacent departure and arrival encounters, denoted by $I_{\mathrm{D},i,j}$ and $I_{\mathrm{A},i,j}$, respectively. Here, D represents departure and A represents arrival. The quantity $v_\infty$, known as the hyperbolic excess velocity, represents the difference between the spacecraft velocity and the velocity of the celestial body. The order in which the module generates legs is arbitrary. However, the computation becomes most efficient when the module first processes the combinations with the largest number of adjacent encounter data.
\begin{equation}
\begin{aligned}
&\mathrm{LD}_{i,j} = [I_{\mathrm{L},i,j}, I_{\mathrm{D},i,j}, I_{\mathrm{A},i,j}, v_{\infty \mathrm{D},i,j}, v_{\infty \mathrm{A},i,j},  H_{\mathrm{leg},i,j}, N_{\mathrm{rev},i,j}],  \\
&\text{Number of legs} :\quad i = 1,\ldots,n_i-1, \\
&\text{Number of data in each leg} :\quad j = 1,\ldots,n_j
\end{aligned}
\end{equation}
In addition, Lambert solutions can become numerically unstable for transfers near $180^\circ$ or $360^\circ$. In such cases, the \textit{Leg Module} first determines the magnitude of the Lambert departure velocity as
\begin{equation}
v_{\mathrm{res}} = | \bm{v}_{\mathrm{D,Lambert}} |.
\end{equation}
The module then keeps this magnitude fixed and analytically computes nearby trajectories to avoid missing feasible solutions. Since infinitely many orbital planes are possible, the departure hyperbolic excess velocity $\bm{v}_{\infty,\mathrm{D}}$ is expressed using two angles, $\kappa$ and $\alpha$, relative to the orbital plane of the departure body as
\begin{equation}
v_\mathrm{D} =
v_\mathrm{res}
\left(
\cos(\alpha)\bm{e_1}
+
\sin(\alpha)\sin(\kappa)\bm{e}_2
+
\sin(\alpha)\cos(\kappa)\bm{e}_3
\right),
\quad
(\bm{e}_1,\bm{e}_2,\bm{e}_3) =
\left(
\frac{\bm{v}_\mathrm{B}}{v_\mathrm{B}},
\frac{\bm{r}_{\mathrm{D}}\times\bm{e}_1}{\|\bm{r}_{\mathrm{D}}\times\bm{e}_1\|},
\bm{e}_1\times\bm{e}_2
\right)
\end{equation}
The angle $\alpha$, which represents the angle between the celestial-body velocity and the Lambert velocity vector, is sampled within the velocity limits defined by the user with step size $\Delta v_\infty$. For full-revolution transfers, namely transfers close to $360^\circ$, the \textit{Leg Module} also treats the crank angle $\kappa$ as a user-defined free parameter.
\begin{equation}
\cos{\alpha}
=
\frac{v_\mathrm{res}^2 + v_\mathrm{B}^2 - v_\infty^2}{2 v_\mathrm{res} v_\mathrm{B}}.
\end{equation}
On the other hand, for transfers near $180^\circ$, referred to as $\pi$ transfers, the module computes the crank angle $\kappa$ so that the radial velocity obtained from the Lambert solution is preserved. The same physical constraint also applies to transfers between different celestial bodies.
\begin{align}
\bm{r}_\mathrm{D}^\intercal \bm{v}_\mathrm{D}
&=
\bm{r}_\mathrm{D}^\intercal \bm{v}_{\mathrm{D,Lambert}}\\
\cos{\kappa}
&=
\frac{
\bm{r}_\mathrm{D}^\intercal
\left(
\bm{r}_\mathrm{D}
-
v_\mathrm{res}\cos(\alpha)\bm{e}_1
\right)
}{
\bm{r}_\mathrm{D}^\intercal\bm{e}_3 \sin(\alpha)
}
\end{align}
The arrival velocity $\bm{v}_\mathrm{A}$ is computed as
\begin{align}
\bm{h} &= \bm{r}_\mathrm{D} \times \bm{v}_\mathrm{D}\\
\bm{v}_\mathrm{A}
&=
\frac{
(\bm{r}_\mathrm{A}^\intercal v_\mathrm{A,Lambert})\bm{r}_\mathrm{A}
+
(\bm{h}\times \bm{r}_\mathrm{A})
}{
\bm{r}_\mathrm{A}^\intercal\bm{r}_\mathrm{A}
}
\end{align}
The remaining position mismatch after this near-resonant construction is corrected using the resonance-targeting maneuver-placement procedure described in the reference ~\cite{Landau2022-fk}. 

During the generation of leg data, three filters are applied to remove infeasible solutions. The indices $p,q$ represent the indices of adjacent encounter data combinations. The first filter removes leg candidates whose time of flight does not satisfy the predefined bounds before solving the Lambert problem.
\begin{equation}
\tau_{\mathrm{min},i,j}
\le
t_{i+1,q} - t_{i,p}
\le
\tau_{\mathrm{max},i,j}, \quad \forall i \in \mathbb{Z}_1^{n_i-1}, \forall j \in \mathbb{Z}_1^{n_j}
\end{equation}
The second filter removes legs whose $v_\infty$ values fall outside the predefined bounds after solving the Lambert problem.
\begin{align}
v_{\infty,\mathrm{min},i}
\le
&v_{\infty\mathrm{D},i,p}
\le
v_{\infty,\mathrm{max},i}, \quad \forall i \in \mathbb{Z}_1^{n_i-1} \\
v_{\infty,\mathrm{min},i+1}
\le
&v_{\infty\mathrm{A},i,q}
\le
v_{\infty,\mathrm{max},i+1}, \quad \forall i \in \mathbb{Z}_1^{n_i-1}
\end{align}
The third filter removes encounter data (ED) entries that are no longer used by any leg after the preceding filters have been applied. In this filter, the algorithm searches for encounter IDs that are not referenced by any remaining leg in each segment and removes the corresponding ED entries. In addition, during the generation of the next segment, the algorithm does not compute legs involving those removed ED entries. Through these operations, the algorithm gradually reduces the number of ED entries that do not contribute to any connection between adjacent segments, thereby decreasing the number of legs that must be computed. Since only a small portion of all possible combinations forms feasible legs in each segment, this pruning process is important for avoiding combinatorial explosion.

\subsubsection{Flyby Module}

The \textit{Flyby Module} connects adjacent leg data to generate flyby data (FD). Each FD entry contains its own ID, the IDs of the incoming and outgoing legs, and the $\Delta v$ associated with the flyby corresponding to the combination of those two legs.
\begin{equation}
\begin{aligned}
&\mathrm{FD}_{i,j} = [I_{\mathrm{F},i,j}, I_{\mathrm{I},i,j}, I_{\mathrm{O},i,j}, \Delta v_{i,j}] , \\
&\text{Number of flybys} :\quad i = 2,\ldots,n_i-1, \\
&\text{Number of data in each flyby} :\quad j = 1,\ldots,n_j
\end{aligned}
\end{equation}

The module examines combinations of adjacent leg data exhaustively. During this process, the module applies filters and generates flyby data only when a combination satisfies the user-defined constraints. The indices $p$ and $q$ denote the indices of the adjacent leg data being combined.

The first filter enforces the time constraint. The two legs correspond to the incoming and outgoing legs with respect to the flyby body. If the endpoints of these legs do not fall within the allowable time range, the module removes the combination. The module obtains the encounter times from the encounter data stored in the leg data. Here, $I_{\mathrm{D},p}=I_{\mathrm{E},i-1}$ and $I_{\mathrm{A},q}=I_{\mathrm{E},i+1}$.
\begin{equation}
\tau_{\mathrm{min},i-1,i+1}
\le
t_{i+1}(I_{\mathrm{E},i+1}) - t_{i}(I_{\mathrm{E},i-1})
\le
\tau_{\mathrm{max},i-1,i+1}, \quad \forall i \in \mathbb{Z}_2^{n_i-1}
\end{equation}
Second, the \textit{Flyby Module} applies the $\Delta v_\mathrm{Mismatch}$ constraint using the flyby geometry
defined in 
Eqs.~\eqref{equ:delta},\eqref{equ:flyby}, and \eqref{equ:beta}.
For flyby $i$, the incoming and outgoing hyperbolic excess velocity
vectors are given by
$\bm{v}_\infty^-=\bm{v}_{\infty\mathrm{A},i-1}$ and
$\bm{v}_\infty^+=\bm{v}_{\infty\mathrm{D},i}$, respectively.
The total $\Delta v$ required during the flyby consists of the velocity change associated with the difference between the incoming and outgoing hyperbolic excess speeds and the additional maneuver required to achieve the remaining bending. Using the law of cosines,
\begin{equation}
\Delta v^2_\mathrm{Mismatch}
= (\|\bm{v}^-_\infty\| - \|\bm{v}^+_\infty\|)^2
+2\|\bm{v}^-_\infty\| \|\bm{v}^+_\infty\| \frac{\sin^2\beta}{1+\cos\beta}
\end{equation}
If this $\Delta v_\mathrm{Mismatch}$ exceeds the allowable value defined by the user, the module removes the corresponding leg combination. Finally, the module removes legs that exceed the allowable range for all flyby combinations. When the module removes such legs, encounter data, leg data, and flyby data that cannot be connected to the current encounter $i$ may appear in both directions. Therefore, the module sequentially removes data that cannot be connected in the incoming direction using the Incoming Leg Filter and also removes data that cannot be connected in the outgoing direction using the Outgoing Leg Filter. Because the \textit{Leg Module} and the \textit{Flyby Module} interact with each other, the framework does not process them independently. Instead, as the \textit{Leg Module} generates leg data, the \textit{Flyby Module} processes them whenever adjacent legs become available.

\subsubsection{Combo Module}

The \textit{Combo Module} constructs complete trajectories by connecting flyby data. Each flyby data entry consists of three encounter data entries: the flyby body and the two encounters before and after it. This structure is referred to as a triplet. A combination of triplets defines one trajectory. Each unit data element $\mathrm{TD}_{i,j}$ stores only the ID of the flyby data. From this ID, the module obtains the corresponding leg data IDs, and from the leg data, it retrieves the three encounter IDs.
\begin{equation}
\begin{aligned}
&\mathrm{TD}_{j} = [\mathrm{TD}_{2,j}, \mathrm{TD}_{3,j}, \ldots, \mathrm{TD}_{i,j}, \mathrm{TD}_{i+1,j}, \ldots, \mathrm{TD}_{n-1,j}],  \\[6pt]
&\text{Number of triplets} :\quad i = 2,\ldots,n_i-1, \\
&\text{Number of data in a trajectory} :\quad j = 1,\ldots,n_j
\end{aligned}
\end{equation}

The process in the \textit{Combo Module} begins by combining triplets at two specific encounter points. During this process, the module applies filters to remove triplets that cannot be connected. Since the number of triplets differs at each encounter, the computation becomes more efficient when the module starts from the encounter pair with the largest number of triplets.

The \textit{Combo Module} applies two filters in this process. The first filter enforces the time constraint. For adjacent triplets, the module estimates the maximum time of flight from the endpoint times $t_{i-1}$ and $t_{i+2}$. If this value exceeds the allowable range defined by the user, the module removes the Triplet combination. The second filter enforces the $\Delta v$ constraint. Each triplet contains the $\Delta v$ value of the flyby, and the module obtains the total $\Delta v$ by summing these values. If the total exceeds the allowable value $\Delta v_{\mathrm{total}}$ defined by the user, the module removes the Triplet combination. Starting from this initial pair, the module sequentially expands the sequence by connecting adjacent triplets. At each step, it prioritizes the direction—either forward or backward—from the current endpoints that offers a larger number of triplet combinations. This sequential connection continues until a continuous, end-to-end TD set spanning all encounter points is successfully generated.
\begin{align}
\Delta v_j =
\Delta v_{i,j}(\mathrm{TD}_{i,j})
+
\Delta v_{i+1,j}(\mathrm{TD}_{i+1,j}), \quad \forall i \in \mathbb{Z}_2^{n_i-1},\forall j \in \mathbb{Z}_1^{n_j}
\end{align}

The algorithm sequentially generates combinations of triplets at adjacent encounters using the previously obtained two-encounter combinations. It prioritizes adjacent encounters with a large number of combinations when determining the processing order. Similarly, the algorithm estimates the maximum TOF from the endpoint times and prunes combinations based on the allowable range defined by the user. At the same time, it accumulates the $\Delta v$ values and also prunes a trajectory if the total $\Delta v$ exceeds the allowable upper limit $\Delta v_{\mathrm{total}}$.
Finally, all obtained $\mathrm{TD}_j$ are converted into the output format $O_j$.
\begin{equation}
\begin{alignedat}{3}
&O_j = [t_{1,j}, E_{1,j}, v_{\infty\mathrm{D},1,j}, v_{\infty\mathrm{A},1,j}, \Delta v_{1,j}, t_{i,j}, E_{i,j}, v_{\infty\mathrm{D},i,j}, v_{\infty\mathrm{A},i,j}, t_{n,j}, E_{n,j}] , \\
&\text{Number of triplets} :\quad i = 2,\ldots,{n_i}-1, \\
&\text{Number of data in a trajectory} :\quad j = 1,\ldots,n_j
\end{alignedat}
\end{equation}

\subsection{Convex Pruning for Broad Search} \label{sec:cp}

Convex programming is an optimization approach for finding an optimal solution to a design problem formulated as a convex problem. Solving such a convex optimization problem enables fast and deterministic computation and guarantees global optimality within the convexified problem. A convex formulation also allows flexible definition of decision variables, objective functions, and constraints. For constrained multivariable problems such as low-thrust transfer trajectory design, previous studies have reported that even a single convex programming iteration can be sufficient to obtain a physically plausible optimized solution when the linearized trajectory propagation is formulated as equality constraints \cite{Ozaki2022-zz}.

In the Broad Search stage of this framework, we introduce convex pruning (CP) into the \textit{Leg Module} of the trajectory generation process.  CP applies convex programming to near-feasible trajectory legs in order to evaluate whether bounded continuous thrust can recover them before the broad-search algorithm discards them. Specifically, when a Lambert solution does not involve a resonant transfer and either of the terminal hyperbolic excess velocities, $v_{\mathrm{D},\infty}$ or $v_{\mathrm{A},\infty}$, exceeds the prescribed allowable range, the corresponding leg is refined using CP. In this pruning, the objective function is defined as the total magnitude of continuous thrust over the trajectory, while the constraints include the linearized equations of motion, the maximum hyperbolic excess velocity at the flyby points, and the upper thrust bound $u_{\max}$. To reduce computational cost, the computation is restricted to cases where the departure and arrival hyperbolic excess velocities, $v_\infty$, satisfy the user-defined boundary ratio $\zeta$. This restriction avoids evaluating a large number of infeasible solutions.

\subsubsection{Overall Computational Procedure}
Although the CP formulation is applicable to arbitrary state
representations, this study adopts MEE and formulates
the equations of motion based on Eqs.~\eqref{equ: dynamical} and
\eqref{equ: dynamical control} as
\begin{equation}
\begin{aligned}
&\dot{\bm{x}}_\mathcal{M} = \bm{f}_{\mathrm{kep},\mathcal{M}}(\bm{x}_\mathcal{M}) + \widetilde{\bm{F}}_\mathcal{M}(\bm{x}_\mathcal{M}) \bm{{u}}_\mathcal{M}\\
&\text{where }\|\bm{u}_{\mathcal{M}}\| \le 1 , \quad \widetilde{\bm{F}}_\mathcal{M}(\bm{x}_\mathcal{M}) = a_{\mathrm{acc},\mathrm{max}} \bm{F}_\mathcal{M}(\bm{x}_\mathcal{M}) \label{equ:mee dynamics}
\end{aligned}
\end{equation}
Here, $a_{\mathrm{acc},\max}$ denotes the fixed upper bound on the perturbing-acceleration magnitude used in the CP formulation and is defined as $a_{\mathrm{acc},\max}
=(D N_{\mathrm{EP}}\,T_{\mathrm{thrust},\mathrm{max}})
/{m_{\mathrm{ref}}}$,
where $T_{\mathrm{thrust},\max}$ is the maximum thrust produced by each unit, and $m_{\mathrm{ref}}$ is a prescribed reference spacecraft mass. We define these quantities before the broad-search phase, and their combination provides a constant acceleration bound for CP.
The state and control variables are defined as 
\begin{align}
\bm{x}_\mathcal{M}=[p,f,g,h,k,L]^\top,
\quad
\bm{u}_\mathcal{M}=[u_R,u_T,u_N]^\top
\end{align}
where $u_R$, $u_T$, and $u_N$ denote the radial, transverse, and
normal components of the dimensionless control vector in the local
RTN frame, respectively. Under unperturbed Keplerian motion, the MEE components
$(p,f,g,h,k)$ remain constant, while only the longitude $L$ evolves \cite{Walker1985-wo}. In addition, MEE avoid the singularities of
classical orbital elements associated with circular and equatorial
orbits~\cite{oguri2022stochastic,Oguri2022-cr}. 
The control influence matrix $\bm{F}_\mathcal{M}(\bm{x}_\mathcal{M})$ and the drift term $\bm{f}_{\mathrm{kep},\mathcal{M}}(\bm{x}_\mathcal{M})$ in the MEE formulation are given by
\begin{equation}
\bm{f}_{\mathrm{kep},\mathcal{M}}(\bm{x}_\mathcal{M}) =
\begin{bmatrix}
0 & 0 & 0 & 0 & 0 & \sqrt{\mu p}\left(\frac{q}{p}\right)^2
\end{bmatrix}^\top
\end{equation}
\begin{equation}
\bm{F}_\mathcal{M} (\bm{x}_\mathcal{M}) =
\begin{bmatrix}
0 & \frac{2p}{q}\sqrt{\frac{p}{\mu}} & 0 \\
\sqrt{\frac{p}{\mu}}\sin L &
\sqrt{\frac{p}{\mu}}\frac{1}{q}[(q+1)\cos L+f] &
-\sqrt{\frac{p}{\mu}}\frac{g}{q}[h\sin L-k\cos L] \\
-\sqrt{\frac{p}{\mu}}\cos L &
\sqrt{\frac{p}{\mu}}\frac{1}{q}[(q+1)\sin L+g] &
\sqrt{\frac{p}{\mu}}\frac{f}{q}[h\sin L-k\cos L] \\
0 & 0 & \sqrt{\frac{p}{\mu}}\frac{s^2\cos L}{2q} \\
0 & 0 & \sqrt{\frac{p}{\mu}}\frac{s^2\sin L}{2q} \\
0 & 0 & \sqrt{\frac{p}{\mu}}\frac{1}{q}[h\sin L-k\cos L]
\end{bmatrix}
\end{equation}
where the auxiliary variables are defined as
\begin{equation}
q = 1 + f\cos L + g\sin L,
\qquad
s^2 = 1 + h^2 + k^2 
\end{equation}
For notational simplicity, the state dependence is suppressed hereafter,
such that
$
\bm{F}_{\mathcal M}(t)
:=
\bm{F}_{\mathcal M}
(\bm{x}_{\mathcal M}(t))$
,
$
\widetilde{\bm{F}}_{\mathcal M}(t)
:=
\widetilde{\bm{F}}_{\mathcal M}
(\bm{x}_{\mathcal M}(t))$, and $
\bm{f}_{\mathrm{kep},\mathcal M}(t)
:=
\bm{f}_{\mathrm{kep},\mathcal M}
(\bm{x}_{\mathcal M}(t)).
$
Linearizing the MEE dynamics in Eq.~\eqref{equ:mee dynamics} about the
reference trajectory
$(\bar{\bm{x}}_{\mathcal M}(t),\bar{\bm{u}}_{\mathcal M}(t))$
yields:
\begin{align}
\dot{\bm{x}}_\mathcal{M}(t) &= A_\mathcal{M}(t) \bm{x}_\mathcal{M}(t) + B_\mathcal{M}(t)\bm{{u}}_\mathcal{M}(t) + \bm{c}_\mathcal{M}(t),\\
A_\mathcal{M}(t) &=
\left.
\frac{\partial (\widetilde{\bm{F}}_\mathcal{M}(t) \bm{u}_\mathcal{M}(t)+\bm{f}_{\mathrm{kep},\mathcal{M}}(t))}{\partial \bm{x}_\mathcal{M}(t)}
\right|_{\bar{\bm{x}}_\mathcal{M}(t),\bar{\bm{u}}_\mathcal{M}}
=
\frac{\partial \widetilde{\bm{F}}_\mathcal{M}(t)}{\partial \bm{x}_\mathcal{M}(t)}\bar{\bm{u}}_\mathcal{M}(t)
+
\frac{\partial \bm{f}_{\mathrm{kep},\mathcal{M}}(t)}{\partial \bm{x}_\mathcal{M}(t)},
\\
B_\mathcal{M}(t) &=
\left.
\frac{\partial (\widetilde{\bm{F}}_\mathcal{M}(t) \bm{u}_\mathcal{M}(t)+\bm{f}_{\mathrm{kep},\mathcal{M}}(t))}{\partial \bm{u}_\mathcal{M}(t)}
\right|_{\bar{\bm{x}}_\mathcal{M}(t),\bar{\bm{u}}_\mathcal{M}(t)}
=
\widetilde{\bm{F}}_\mathcal{M}(t),
\\
\bm{c}_\mathcal{M}(t) &= \bm{f}_{\mathrm{kep},\mathcal{M}}(t)(\bar{\bm{x}}_\mathcal{M}(t),\bar{\bm{u}}_\mathcal{M}(t))
- A_\mathcal{M}(t)\bar{\bm{x}}_\mathcal{M}(t)
- B_\mathcal{M}(t)\bar{\bm{u}}_\mathcal{M}(t)
\end{align}
Following \cite{Aziz2019-ob}, the state transition matrix (STM) is approximated by multiplying the instantaneous sensitivity by the time step, $\tau$. 
\begin{equation}
\Phi_\mathcal{M}(t_{j+1},t_j) = \frac{\delta \bm{x}_{\mathcal{M}}(t_{j+1})}{\delta \bm{x}_{\mathcal{M}}(t_j)}
\approx
I_{6 \times 6} + A_\mathcal{M}(t_j) \tau_j, \quad \tau_j = t_{j+1}-t_j, \quad \forall j \in \mathbb{Z}_1^{N_j-1}
\end{equation}
Here, this study uses the zero-order hold (ZOH) approximation, in which the control components are assumed to be constant within each segment. $N_j$ denotes the number of discretization nodes along the leg.
Using these approximations, the linearized dynamics are discretized over each time segment $[t_j,t_{j+1}]$ as
\begin{equation}
\begin{aligned}
\bm{x}_{\mathcal{M},j+1} &= A_{\mathcal{M},j} \bm{x}_{\mathcal{M},j} + B_{\mathcal{M},j} \bm{u}_{\mathcal{M},j} + \bm{c}_{\mathcal{M},j} , \quad \forall j \in \mathbb{Z}_1^{N_j-1} \label{equ:discrete-dynamics}
 \\
A_{\mathcal{M},j} &= \Phi_\mathcal{M}(t_{j+1},t_j)
\approx I_{6 \times 6} +
\left(
\frac{\partial \widetilde{\bm{F}}_{\mathcal{M}}(t_j)}{\partial \bm{x}_\mathcal{M}(t_j)}\bar{\bm{u}}_{\mathcal{M},j}
+
\frac{\partial \bm{f}_{\mathrm{kep},\mathcal{M}}(t_j)}{\partial \bm{x}_\mathcal{M}(t_j)}
\right)\tau_j\\
B_{\mathcal{M},j} &= \int_{t_j}^{t_{j+1}} \Phi_\mathcal{M}(t_{j+1},t) B_\mathcal{M}(t)  \, dt  \approx \widetilde{\bm{F}}_{\mathcal{M}}(t_j)\tau_j, \\
\bm{c}_{\mathcal{M},j} &= \bar{\bm{x}}_{\mathcal{M},j+1}-A_{\mathcal{M},j}\bar{\bm{x}}_{\mathcal{M},j}-B_{\mathcal{M},j}\bar{\bm{u}}_{\mathcal{M},j}
\end{aligned}
\end{equation}
At each encounter point, the spacecraft boundary state can be expressed
by adding the spacecraft hyperbolic excess velocity vector to the
velocity component of the encounter-body state. However, in the MEE formulation, the position and velocity components are not explicitly separated. Here, the instantaneous sensitivity matrix makes it possible to incorporate the velocity components separately. At the boundary points, the formulation uses the same matrix to map an
instantaneous RTN velocity increment, rather than an acceleration, to the
corresponding variation in the MEE state. The boundary conditions are written as
\begin{align}
\bm{x}_{\mathcal{M},1} &= \bm{x}_{\mathcal{M},E}(t_\mathrm{D}) +  \|\bm{v}_{\infty\mathrm{D},\max}\|\bm{F}_{\mathcal{M}}(t_\mathrm{D})\bm{u}_{\mathcal{M},\mathrm{D}} \label{equ:con_x_0}\\ 
\bm{x}_{\mathcal{M},N_J} &= \bm{x}_{\mathcal{M},E}(t_\mathrm{A}) +  \|\bm{v}_{\infty\mathrm{A},\max}\|{\bm{F}}_{\mathcal{M}}(t_\mathrm{A})\bm{u}_{\mathcal{M},\mathrm{A}} \label{equ:con_x_f}
\end{align}
Here, the subscripts $\mathrm{D}$ and $\mathrm{A}$ denote the departure
and arrival encounter points of the trajectory leg, respectively.
The vector $\bm{x}_{\mathcal M,E}$ denotes the MEE state of the
encounter body. 
$\|\bm{v}_{\infty\mathrm{D},\max}\|$ and $\|\bm{v}_{\infty\mathrm{A},\max}\|$ denote the user-specified upper bounds on the hyperbolic excess velocities at the departure and arrival encounters, respectively. The products
$\|\bm{v}_{\infty\mathrm{D},\max}\|\bm{u}_{\mathcal M,\mathrm{D}}$ and
$\|\bm{v}_{\infty\mathrm{A},\max}\| \bm{u}_{\mathcal M,\mathrm{A}}$ represent
the allowable departure and arrival hyperbolic excess velocity vectors,
respectively, expressed in the local RTN frame. The matrix
$\bm{F}_{\mathcal M}(t)$ provides the instantaneous first-order mapping
from the RTN velocity increment to the corresponding variation in the
MEE state.
In addition, the formulation constrains the thrust-control vector and the boundary
hyperbolic-excess-velocity direction vectors as
\begin{equation}
\|\bm{u}_{\mathcal{M},\mathrm{D}}\| \le 1,
\qquad
\|\bm{u}_{\mathcal{M},\mathrm{A}}\| \le 1,
\qquad
\|\bm{u}_{\mathcal{M},j}\| \le 1, \quad \forall j \in \mathbb{Z}_1^{N_j-1}
\label{equ:con_u_max}
\end{equation}
Finally, the CP formulates the convex optimization problem by combining the objective function and all constraints:
\begin{equation}
\begin{aligned} 
&\min_{\substack{ \bm{x}_{\mathcal{M},j},{\forall j \in \mathbb{Z}_1^{N_j}} , \bm{u}_{\mathcal{M},j},{\forall j \in \mathbb{Z}_1^{N_j-1}},\\\bm{u}_{\mathcal{M},\mathrm{D}} , \bm{u}_{\mathcal{M},\mathrm{A}}} } && \sum_{j=1}^{N_j-1}  
\|\bm{u}_{\mathcal{M},j}\| \tau_j  \\
&\text{subject to} &&\bm{x}_{\mathcal{M},j+1} = A_{\mathcal{M},j} \bm{x}_{\mathcal{M},j} + B_{\mathcal{M},j} \bm{u}_{\mathcal{M},j} + \bm{c}_{\mathcal{M},j}, \quad   \forall j \in \mathbb{Z}_1^{N_j-1}\\
& &&\bm{x}_{\mathcal{M},1} = \bm{x}_{\mathcal{M},E}(t_\mathrm{D}) +  \|\bm{v}_{\infty\mathrm{D},\max}\|\bm{F}_{\mathcal{M}}(t_\mathrm{D})\bm{u}_{\mathcal{M},\mathrm{D}} \\
& &&\bm{x}_{\mathcal{M},N_j} = \bm{x}_{\mathcal{M},E}(t_\mathrm{A}) +  \|\bm{v}_{\infty\mathrm{A},\max}\|{\bm{F}}_{\mathcal{M}}(t_\mathrm{A})\bm{u}_{\mathcal{M},\mathrm{A}}  \\
& &&\|\bm{u}_{\mathcal{M},\mathrm{D}}\| \le 1, \quad \|\bm{u}_{\mathcal{M},\mathrm{A}}\| \le 1,
\quad \|\bm{u}_{\mathcal{M},j}\| \le 1 ,\quad  \forall j \in \mathbb{Z}_1^{N_j-1}  
\end{aligned}
\label{equ:cp optimization}
\end{equation}

where $t_\mathrm{D}$ and $t_\mathrm{A}$ denote the fixed departure and
arrival epochs, respectively. The quantities
$\bm{x}_{\mathcal{M},E}(t_\mathrm{D})$,
$\bm{x}_{\mathcal{M},E}(t_\mathrm{A})$,
$\|\bm{v}_{\infty\mathrm{D},\max}\|$,
$\|\bm{v}_{\infty\mathrm{A},\max}\|$,
$\bm{F}_{\mathcal{M}}(t_\mathrm{D})$, and
$\bm{F}_{\mathcal{M}}(t_\mathrm{A})$
are fixed parameters in the convex subproblem. The linearized dynamics matrices
$A_{\mathcal{M},j}$, $B_{\mathcal{M},j}$, and $\bm{c}_{\mathcal{M},j}$
are also constant for each $j$ during each convex-programming solve.

\subsubsection{Convex Pruning Accuracy}
To evaluate the accuracy of the proposed CP formulation, this study compares the CP-based trajectory against a trajectory propagated using the solution obtained by CP via conventional ordinary differential equation (ODE) integration.
As an example, one trajectory leg is selected for detailed analysis from the initial solution associated with the minimum-time-of-flight representative trajectory presented in Section~\ref{sec:representative}, namely the solution obtained with CP and labeled ID~10. This leg corresponds to the transfer from event 3, a Venus flyby, to event 4, another Venus flyby, in the mission sequence. For this analysis, we discretize the selected trajectory leg into $N_j=31$ nodes. Additionally, for this calculation, Eq.~\eqref{equ:cp optimization} sets the hyperbolic excess velocity bounds at both the departure and arrival encounters to $\| \Delta \bm{v}_{\infty \mathrm{D},\mathrm{max}}
\| = \| \Delta \bm{v}_{\infty \mathrm{A},\mathrm{max}}
\|=9 \, \mathrm{[km/s]}$.

\begin{figure}[!htbp]
\centering
\begin{minipage}{0.48\textwidth}
\centering
\includegraphics[width=\linewidth]{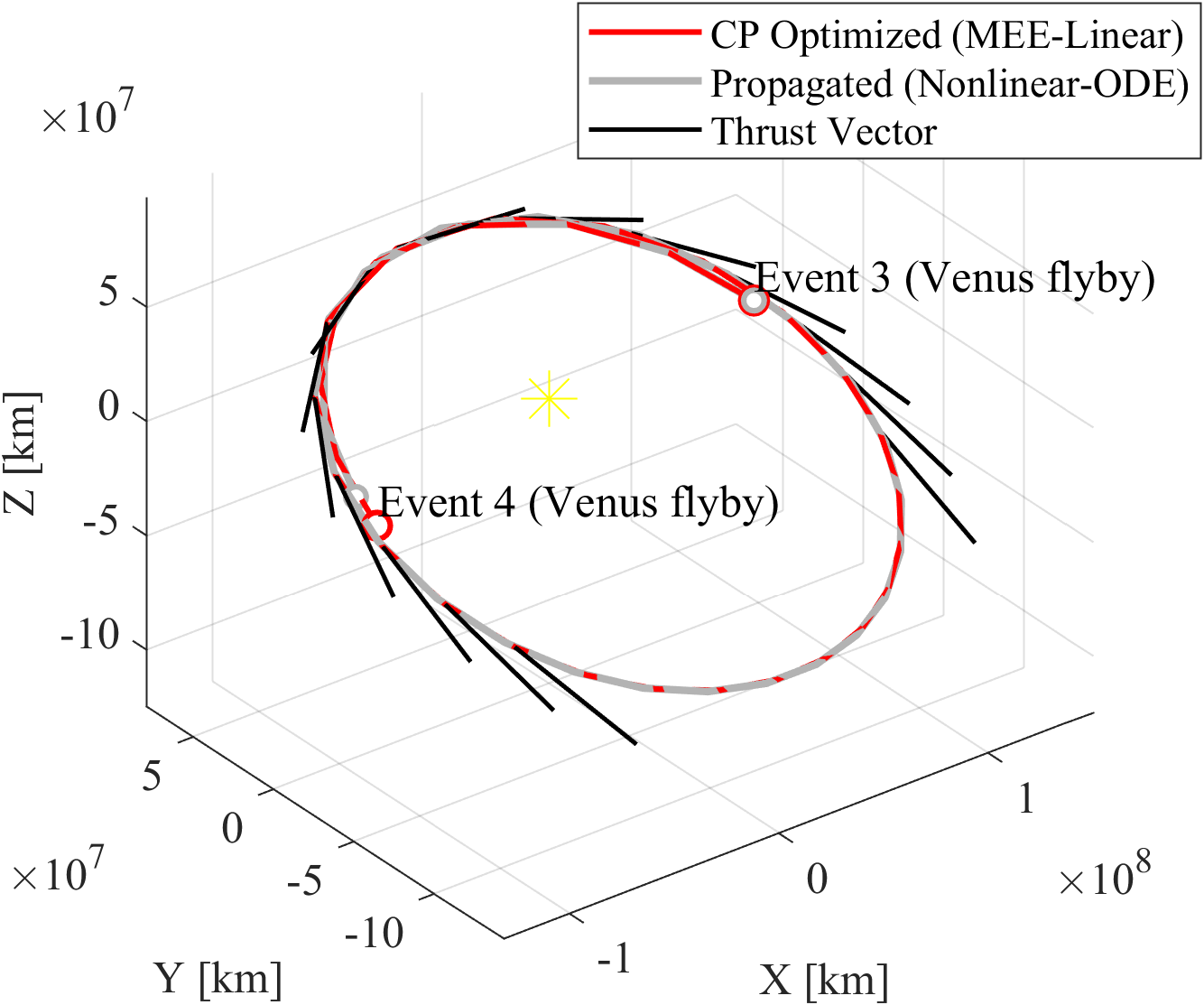}
\end{minipage}\hfill
\begin{minipage}{0.48\textwidth}
\centering
\includegraphics[width=\linewidth]{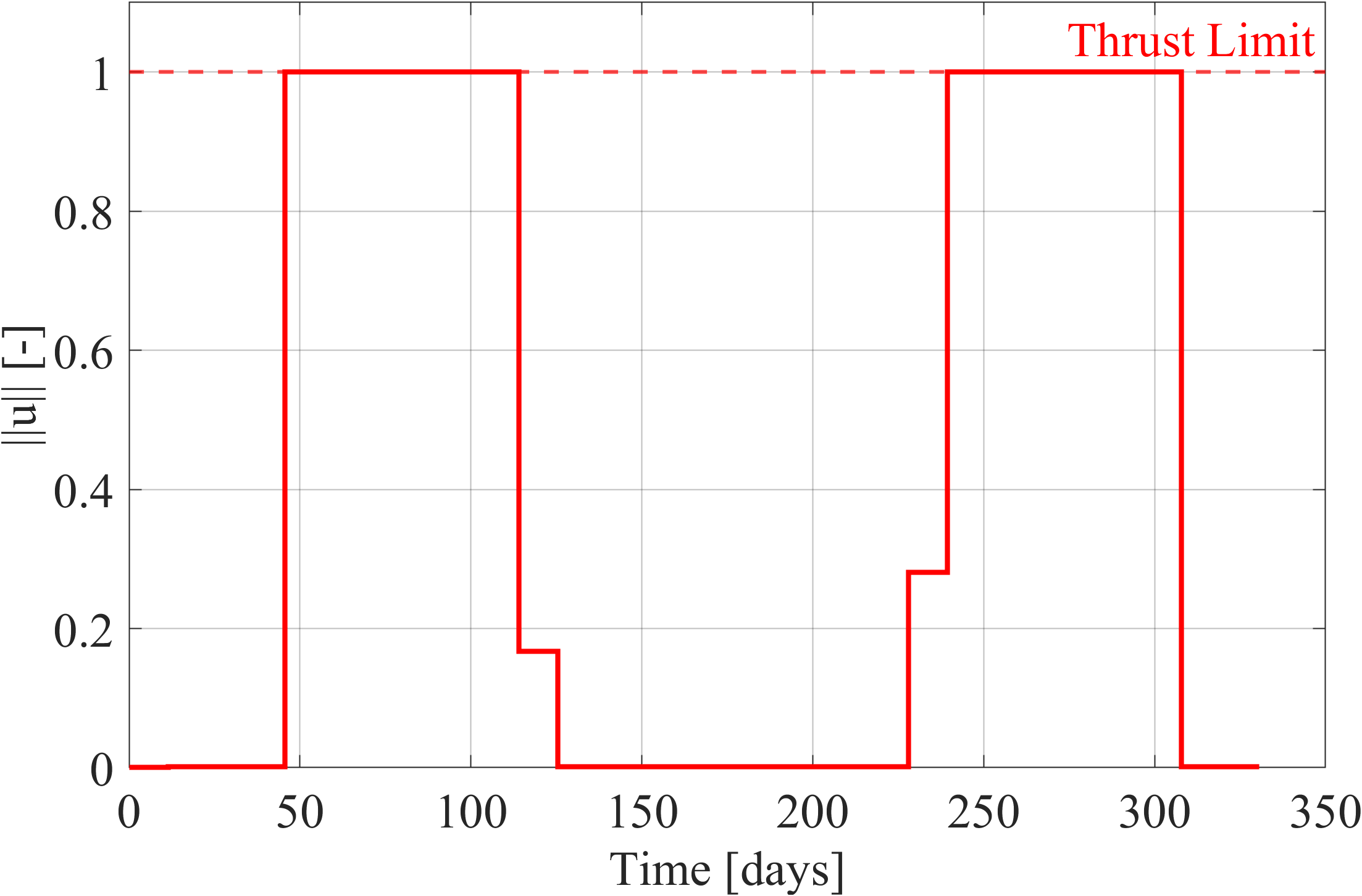}
\end{minipage}

\caption{
(Left) Comparison of the CP-optimized trajectory and the nonlinear ODE propagation. (Right) The corresponding optimized thrust profile.
}
\label{fig:CPODE_Comparison}
\end{figure}

Fig.~\ref{fig:CPODE_Comparison} presents a three-dimensional comparison between the CP solution and the nonlinearly propagated trajectory, together with the associated low-thrust control profile.
As shown in Fig.~\ref{fig:CPODE_Comparison}, the optimized control profile exhibits a predominantly bang--bang-like structure, with the thrust acceleration operating mainly near either zero or its prescribed upper limit. This behavior is consistent with the thrust-saturation structure expected from Pontryagin's maximum principle for a minimum-fuel low-thrust trajectory problem. The result suggests that the convex program captures a physically meaningful control structure that approaches the theoretical structure of an optimal low-thrust solution. Additionally, CP satisfies the encounter constraints of the hyperbolic excess velocity bounds by introducing approximately 0.8 km/s of total low-thrust. Before applying CP, the Lambert solution yields a hyperbolic excess velocity of $12.3 \, \mathrm{km/s}$ at the first encounter, which exceeds the prescribed upper bound value of $9.0 \, \mathrm{km/s}$ by $3.3 \, \mathrm{km/s}$. The CP result demonstrates that the trajectory can satisfy the encounter-velocity constraint with a substantially smaller accumulated low-thrust than this initial excess-velocity difference.

\begin{figure}
    \centering
    \includegraphics[width=1.0\linewidth]{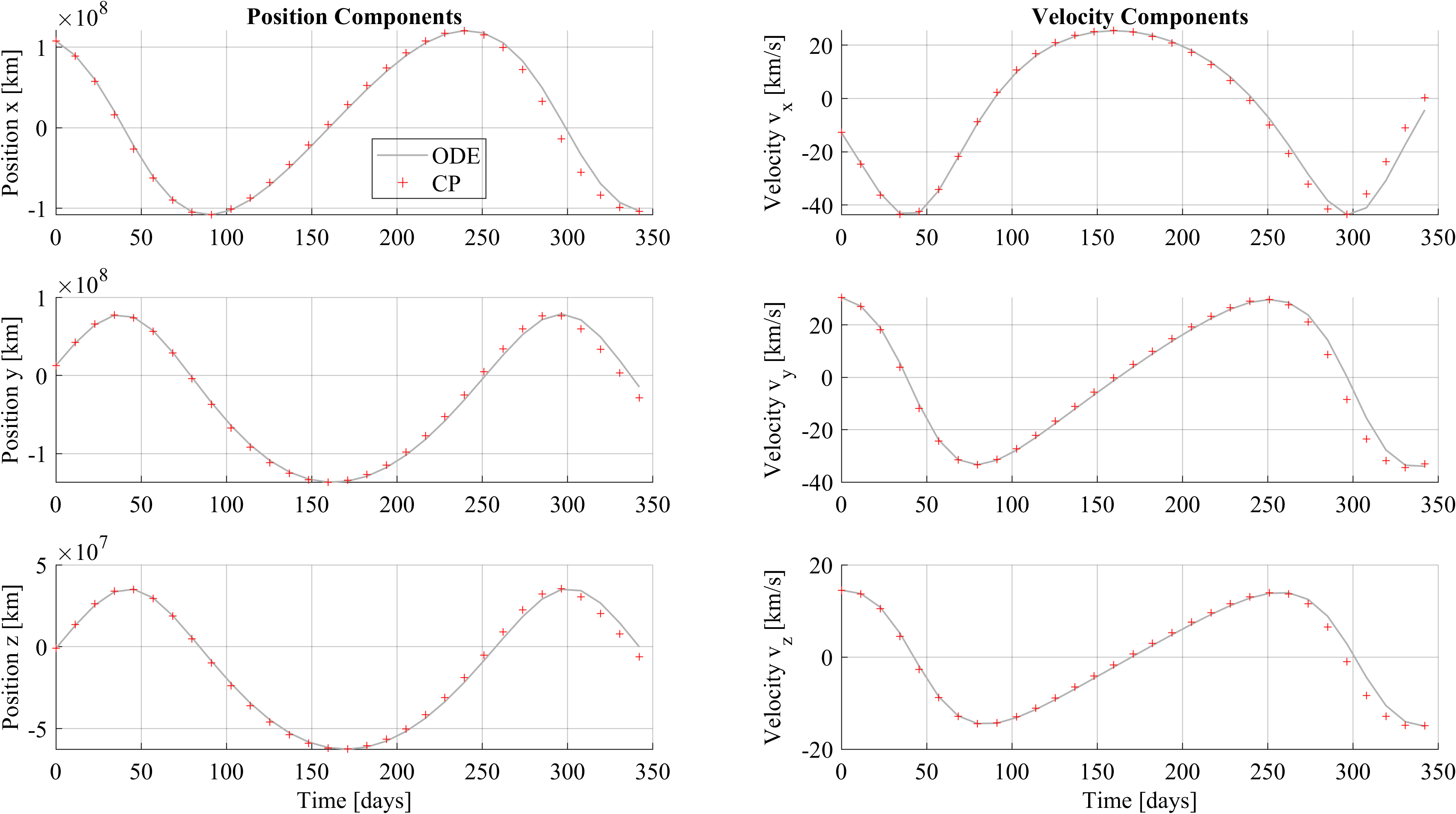}
    \caption{Comparison of state vector components between the CP solution and ODE-based propagation for a representative trajectory leg in the Broad Search.}
    \label{fig:StateVector_CP}
\end{figure}

The discrepancy between the CP solution and the trajectory obtained through nonlinear ODE propagation depends on the magnitude of the correction required for CP to reduce the Lambert-solution hyperbolic excess velocity to the prescribed bounds, $\|\bm{v}_{\infty\mathrm{D},\max}\|$ and $\|\bm{v}_{\infty\mathrm{A},\max}\|$. Fig.~\ref{fig:StateVector_CP} compares the position and velocity state components for the trajectory. In this relatively demanding recovery case, the propagation discrepancy reaches approximately 10\% of the characteristic trajectory scale. By contrast, when CP requires a correction of approximately $2.0 \, \mathrm{km/s}$ or less to bring the hyperbolic excess velocity of the Lambert solution within the prescribed encounter bounds, the discrepancy between the CP solution and nonlinear ODE propagation remains below 1\%. These results indicate that the approximation accuracy of the CP solution depends strongly on the magnitude of the required velocity correction. Therefore, the single-iteration CP solution is not intended to serve directly as a high-fidelity propagated trajectory; rather, it identifies promising regions of the broad-search space and provides initial estimates for the subsequent nonlinear programming (NLP) stage. Even for the case shown in Fig.~\ref{fig:StateVector_CP}, which exhibits a propagation discrepancy of approximately $10\%$, the subsequent NLP process eliminates the state mismatch and converges to a feasible solution while largely preserving the structure of the low-thrust control profile generated by the CP stage.
These results indicate that the proposed CP approach achieves sufficient dynamical consistency for broad-search pruning while generating physically reasonable control structures at a substantially lower computational cost than a full nonlinear programming solve.

 The proposed method is implemented in MATLAB. In its baseline implementation, the method requires approximately 0.3 sec per solve. Most of the computational cost comes from constructing the linearized dynamical constraints, whereas the optimization itself converges in approximately 0.2 sec. The implementation further accelerates the computation by applying code generation (MEX) to the linearized dynamical constraints and by compiling the YALMIP setup into C programming language using the \texttt{optimizer} interface. With these improvements, the total computation time decreases to approximately 0.05 sec. In the current implementation, the proposed CP typically required $O(10^{-2})$ sec per solve, depending on the transfer geometry and solver behavior.

When the constraints are overly restrictive and no feasible solution exists, the solver terminates rapidly, allowing the search process to proceed without stagnation. This property is particularly advantageous when the method is integrated into a broad search framework such as \textit{Star}, because the framework can efficiently discard infeasible candidates and immediately proceed to the next candidate solution. 

\subsection{Local Optimization with Thrust Regularization} \label{sec:local optimization}

\subsubsection{Overall Optimization Procedure}

The Sims--Flanagan transcription (SFT), introduced by Sims and Flanagan in 1997, is a direct transcription method for low-thrust gravity-assist trajectory optimization \cite{sims1997preliminary}. Because SFT replaces continuous low-thrust propagation with a sequence of Keplerian arcs connected by impulsive velocity increments, it avoids the numerical integration of the low-thrust equations of motion and therefore offers relatively low computational cost. This feature has led to its widespread use in preliminary trajectory design and optimization \cite{Englander2017-ls}.

SFT divides each trajectory leg into a prescribed number of segments. Within each segment, the spacecraft follows an unperturbed Keplerian arc, and an impulsive velocity increment is applied at the segment midpoint to approximate the effect of continuous thrust. Forward propagation from the departure boundary and backward propagation from the arrival boundary generate two half-trajectories, which are matched at an interior patch point.
In conventional SFT formulations, the Cartesian components of $\Delta\bm{v}_\mathcal{C}$ serve directly as optimization variables. For problems whose optimal thrust profile contains coast arcs, however, this parameterization can introduce a singularity into the objective-function and constraint gradient, $\partial \| \Delta \bm{v}_\mathcal{C} \| /\partial \Delta \bm{v}_\mathcal{C} = \Delta \bm{v}_\mathcal{C}^\top / \| \Delta \bm{v}_\mathcal{C}\| $, where the denominator approaches zero for coast segments as the optimization converges toward the optimal solution \cite{Oshima2023-gg,Ottesen2021-cr,Varghese2026-qt}.

To address this issue, this study introduces a thrust-regularization method inspired by~\cite{Oshima2023-gg} that removes the singularity analytically without requiring an artificial tuning parameter. In this method, the velocity increment is reformulated using a variable $\bm{U}$, and the formulation defines $\Delta  \bm{v}_\mathcal{C}$ so that its norm is equal to the squared norm of $\bm{U}$. The formulation of this thrust regularization is as follows:
\begin{equation}
\begin{aligned}
&\bm{U} :=
\begin{bmatrix}
U_x \\ U_y \\ U_z
\end{bmatrix},
\qquad
\Delta \bm{v}_\mathcal{C} :=
\begin{bmatrix}
\Delta v_x \\ \Delta v_y \\ \Delta v_z
\end{bmatrix}
=
\begin{bmatrix}
U_x^2-U_y^2-U_z^2 \\
2 U_x U_y \\
2 U_x U_z
\end{bmatrix},
\qquad
\|\Delta \bm{v}_\mathcal{C}\| = \|\bm{U}\|^2, \\
&\frac{\partial \| \Delta \bm{v}_\mathcal{C} \|} {\partial \bm{U}} = \frac{\partial \| \bm{U} \|^2} {\partial \bm{U}} = 2\bm{U}^\top
\label{dv_reg}
\end{aligned}
\end{equation}
Thus, the gradient no longer contains the velocity-increment norm in the denominator and remains well defined at $U=0$, thereby removing the singularity associated with coast segments.

\subsubsection{Interfacing CP Results with Thrust-Regularized SFT}
Before local optimization, a refinement procedure maps the CP-derived control history onto the SFT discretization and improves its consistency with the nonlinear propagation model. The procedure first resamples the CP-derived control history when the CP and SFT discretizations differ, assigning the control corresponding to the same epoch to each SFT segment. The departure state obtained from CP is then transformed from MEE to Cartesian coordinates and propagated forward using the Keplerian propagation employed in SFT. At the midpoint of each segment, the CP-derived control is transformed into the Cartesian frame using the propagated state and applied as the corresponding velocity increment. The resulting Cartesian state and control histories then define an updated reference trajectory for a second CP solve on the SFT discretization. For this solve, the propagated Cartesian states and control vectors are transformed back into the MEE and RTN representations, respectively, and used as the reference histories $\bar{\bm{x}}_\mathcal{M}$ and $\bar{\bm{u}}_\mathcal{M}$. Section~\ref{sec:cp} discusses the resulting improvement in the CP. After the second CP solve, the refined low-thrust control history is propagated once more using the same Keplerian propagation procedure to obtain the Cartesian velocity increments $\Delta v_\mathcal{C}$ at the SFT segment midpoints. The resulting Cartesian velocity increments provide the control initial guess for local optimization and are converted into the intermediate parameter vector $\bm{U}$. Expressing the Cartesian velocity increment in spherical coordinates gives
\begin{equation}
\begin{bmatrix}
\Delta v_x \\
\Delta v_y \\
\Delta v_z
\end{bmatrix}
=
\|\Delta  \bm{v}_\mathcal{C}\|
\begin{bmatrix}
\cos\psi \\
\sin\psi \cos\phi \\
\sin\psi \sin\phi
\end{bmatrix}, \quad
\psi = \cos^{-1}\left( \frac{\Delta v_x}{\|\Delta  \bm{v}_\mathcal{C}\|} \right), \quad
\phi = \tan^{-1}\left( \frac{\Delta v_z}{\Delta v_y} \right)
\end{equation}
Here, the angular coordinate $\phi$ is evaluated using the two-argument arctangent function, $\mathrm{atan2}(Y,X)$, which avoids the undefined form associated with $\tan^{-1}(0/0)$, with $\tan^{-1}(0,0)=0$ following the MATLAB convention. Even when both $\Delta \bm{v}_y=0$ and $\Delta \bm{v}_z=0$, this implementation returns a finite value for $\phi$ and prevents numerical failure. In addition, the implementation includes a fail-safe for zero-length vectors by replacing a zero norm with a small positive value before evaluating $\psi$ with the inverse cosine. This treatment prevents division by zero in the computation of $\psi$ and further improves numerical robustness.
This representation gives the parameters $U_x$, $U_y$, and $U_z$ as
\begin{align}
\begin{bmatrix}
U_x \\ U_y \\ U_z 
\end{bmatrix}
=
\sqrt{\|\Delta \bm{v}_\mathcal{C} \|}
\begin{bmatrix}
 \cos\frac{\psi}{2} \\
 \sin\frac{\psi}{2}\cos\phi \\
\sin\frac{\psi}{2}\sin\phi
\end{bmatrix} \label{dv_reg}
\end{align}
 
\subsubsection{Thrust-Regularized SFT Formulation}
During local optimization, the SFT divides each trajectory leg into a prescribed number of segments. For leg $i$, the method propagates the trajectory forward from the departure boundary and backward from the arrival boundary through analytical two-body Keplerian arcs. An impulsive velocity increment is applied at the midpoint of each segment, and the forward and backward trajectories are matched at an interior patch point. 
For segment $j$ of trajectory leg $i$, 
the transcription divides each leg into uniformly spaced segments. Therefore, the time interval between consecutive impulse epochs remains constant within leg $i$:
$t_{i,j+1}^{-}-t_{i,j}^{+}=
\tau_{i,j}=(t_{i,N_j}-t_{i,1})/(N_j - 1)$,
where $N_j$ denotes the number of nodes in each leg. Because SFT models each thrust maneuver as an instantaneous impulse, the epochs immediately before and after the impulse coincide:
$t_{i,j}^{-}=t_{i,j}^{+}=t_{i,j}.$
The forward and backward propagation based on the Cartesian formulation introduced in Eq.~\ref{equ:kepler} can be expressed as
\begin{equation}
\begin{aligned}
\text{Forward: }
\begin{cases}
\bm{x}_{\mathcal{C},i,j}^{-}  = \mathcal{K}_\mathcal{C}\left( \bm{x}_{\mathcal{C},i,j-1}^{+},  {\frac{1}{2}}\tau_{i,j} \right), \\ 
\bm{x}_{\mathcal{C},i,j}^{+} = \bm{x}_{\mathcal{C},i,j}^{-}
+ \begin{bmatrix} 
\bm{0}_{3\times 1} \\ \Delta\bm{v}_{\mathcal{C},i,j}(\bm{U}_{i,j}) 
\end{bmatrix},\quad  m_{i,j}^{+} = m_{i,j}^{-} \exp\left(-\frac{\|\bm{U}_{i,j}\|^2}{g_0 I_{\mathrm{sp}}(t_{i,j})}\right),\quad \forall i \in \mathbb{Z}_1^{N_i},\forall j \in \mathbb{Z}_1^{N_j-1}\\ 
\bm{x}_{\mathcal{C},i,j+1}^{-} = \mathcal{K}_\mathcal{C}\left( \bm{x}_{\mathcal{C},i,j}^{+}, {\frac{1}{2}}\tau_{i,j} \right)
\end{cases}
\end{aligned}
\label{equ:sft forward}
\end{equation} 
\begin{equation}
\begin{aligned}
\text{Backward: }
\begin{cases}
\bm{x}_{\mathcal{C},i,j}^{+}  = \mathcal{K}_\mathcal{C}\left( \bm{x}_{\mathcal{C},i,j+1}^{-}, -{\frac{1}{2}}\tau_{i,j} \right), \\ 
\bm{x}_{\mathcal{C},i,j}^{-} = \bm{x}_{\mathcal{C},i,j}^{+}
- \begin{bmatrix} 
\bm{0}_{3\times 1} \\ \Delta\bm{v}_{\mathcal{C},i,j}(\bm{U}_{i,j}) 
\end{bmatrix},\quad  m_{i,j}^{-} = m_{i,j}^{+} \exp\left(\frac{\|\bm{U}_{i,j}\|^2}{g_0 I_{\mathrm{sp}}(t_{i,j})}\right),\quad \forall i \in \mathbb{Z}_1^{N_i},\forall j \in \mathbb{Z}_1^{N_j-1}\\ 
\bm{x}_{\mathcal{C},i,j-1}^{+} = \mathcal{K}_\mathcal{C}\left( \bm{x}_{\mathcal{C},i,j}^{-}, -{\frac{1}{2}}\tau_{i,j} \right)
\end{cases}
\end{aligned}
\label{equ:sft backward}
\end{equation}
where $\mathcal{K}_\mathcal{C}(\bm{x}_\mathcal{C},\pm{\frac{1}{2}}\tau_{i,j})$ denotes the analytical two-body Keplerian propagation over the time interval $\pm{\frac{1}{2}}\tau_{i,j}$ for forward and backward propagation, under the gravitational parameter $\mu$, using a universal-variable formulation
based on Stumpff functions \cite{Shepperd1985-om}. 
The impulsive maneuver updates the spacecraft mass according to the Tsiolkovsky rocket equation. Based on the maximum low-thrust acceleration magnitude, $\|\bm{a}_\mathrm{acc}\|=D N_\mathrm{EP} T_\mathrm{thrust}
(t)/m$, given by Eq.~\eqref{equ: dynamical control}, the thrust-regularization formulation constrains the segmentwise velocity increment during forward and backward propagation using the control variable $\bm{U}_{i,j}$ as
\begin{equation}
\text{Forward: }
\|\bm{U}_{i,j}\|^2 \le \frac{D N_\mathrm{EP} T_\mathrm{thrust} (t_{i,j}) \tau_{i,j}}{m_{i,j}^{+}} ,\quad \forall i \in \mathbb{Z}_1^{N_i},\forall j \in \mathbb{Z}_1^{N_j-1} \label{equ:forward thrust constraint}
\end{equation}
\begin{equation}
\text{Backward: }
\|\bm{U}_{i,j}\|^2 \le \frac{D N_\mathrm{EP}  T_\mathrm{thrust} (t_{i,j}) \tau_{i,j}}{m_{i,j}^{-}}, \quad  \forall i \in \mathbb{Z}_1^{N_i},\forall j \in \mathbb{Z}_1^{N_j-1} \label{equ:backward thrust constraint}
\end{equation}
In addition, the midpoint defect constraints for position, velocity, and mass are defined as follows: 
\begin{equation}
\begin{bmatrix}
\bm{x}_{\mathcal{C},i,\mathrm{mid}}^{\mathrm{+}}-\bm{x}_{\mathcal{C},i,\mathrm{mid}}^{\mathrm{-}}\\
m_{i,\mathrm{mid}}^{\mathrm{+}}-m_{i,\mathrm{mid}}^{\mathrm{-}}
\end{bmatrix}
=\bm{0}_{7\times 1}, \quad \forall i \in \mathbb{Z}_1^{N_i}
\label{equ:sft midpoint constraint}
\end{equation}
where "mid" denotes the common interior patch point of leg $i$, located at the interface between the forward and backward-propagated branches.

The local optimization formulates the trajectory-design problem by combining the SFT with thrust regularization and the gravity-assist constraints introduced in Section~\ref{sec:gravity assist}. The resulting nonlinear programming problem is defined as follows to maximize the final spacecraft mass while enforcing the forward and backward-propagation equations, impulsive velocity and mass updates, center-patch continuity conditions, gravity-assist feasibility conditions, time-of-flight bounds, and thrust constraints:
\begin{equation}
\begin{aligned}
&\max_{\bm{x}_\mathrm{opt}} \quad
&&m_{N_{i},N_j} =m_{1,1} - \sum_{i=1}^{N_i}(m_{i,1}-m
_{i,N_j}) \\
&\mathrm{subject \: to}
&&\text{Dynamics: } \text{Eqs.~\eqref{equ:sft forward} and ~\eqref{equ:sft backward} for forward and backward propagation}, \\ 
& &&\text{Thrust constraints: } 
\text{Eqs.~\eqref{equ:forward thrust constraint} and ~\eqref{equ:backward thrust constraint} for forward and backward propagation}, \\
& &&\text{Patch-point continuity constraints: } 
\text{Eq.~\eqref{equ:sft midpoint constraint}},\\
& &&\text{Gravity-assist constraints: } \text{Eqs.~\eqref{equ:flyby-constraint-regularized} and ~\eqref{equ:velocity manitude}} ,\\
& &&\text{Launch constraints: } m_{1,1} \leq m_\text{delivered}({C_3}_{1,1}), \quad {C_3}_{1,1} \le {C_3}_\mathrm{max}, \\
& &&\text{Time-of-flight constraint: }   \sum_{i=1}^{N_i} (t_{i+1} - t_i) \leq \tau_\mathrm{total,max},\\
&\mathrm{where}
&&\bm{x}_\mathrm{opt}=\{t_{i,1},\bm{v}_{\infty,i,1}^+,m_{i,1}^+,\bm{v}_{\infty,i,N_j}^-,m_{i,N_j}^-
,t_{N_i,N_j},\bm{U}_{i,j}\}, \\
& && \text{Number of Legs: }  i=1,2,...,N_i  \\
& && \text{Number of Nodes for Each Leg: }  j=1,2,...,N_j 
\end{aligned}
\end{equation}
Here, ${C_3}_\mathrm{max}$ denotes the maximum launch characteristic energy, while $\tau_\mathrm{total,max}$ represents the maximum time-of-flight.

\subsubsection{Thrust-Regularization Validation}
\begin{figure} [!t]
    \centering
    \includegraphics[width=0.8\linewidth]{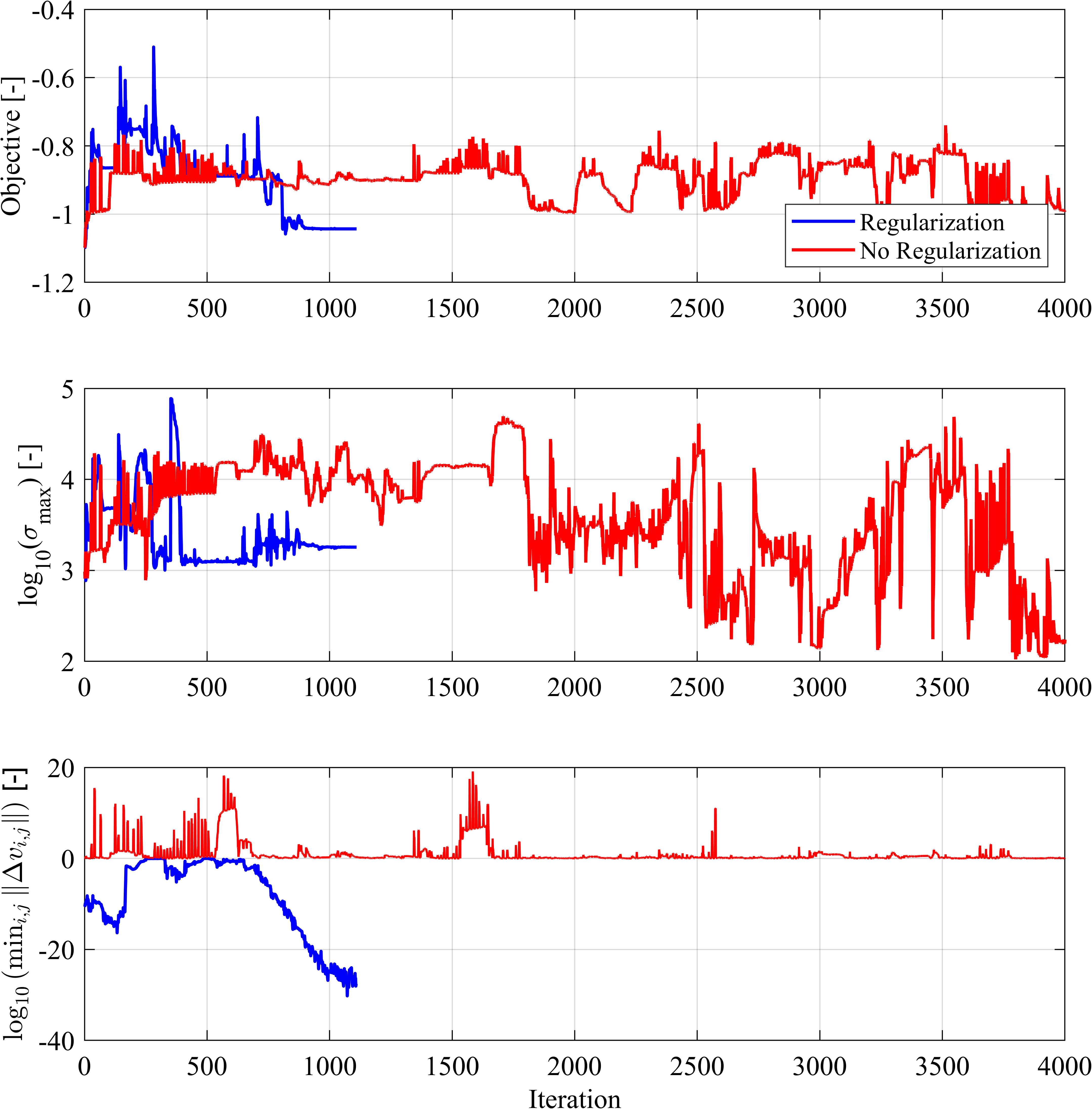}
    \caption{Comparison of optimization convergence characteristics with and without thrust regularization for the test case, corresponding to the w/ CP solution with ID 12, described in Section \ref{sec:example}.} 
    \label{fig:IterationCNRegNoReg}
\end{figure}

To evaluate the effect of introducing thrust regularization into the SFT, we compare the proposed formulation with a conventional SFT formulation that does not employ thrust regularization. For detailed analysis, we select the maximum-final-mass representative trajectory presented in Section~\ref{sec:example}, which corresponds to the w/ CP solution with ID 12.
Fig.~\ref{fig:IterationCNRegNoReg} compares the iteration histories of the optimizations with and without thrust regularization. The plotted quantities include the objective function value, the spectral norm of the full constraint Jacobian, $\sigma_{\max}$, and the minimum control magnitude over the trajectory at each NLP iteration.

The two formulations show clear differences. With thrust regularization, the objective function value largely converges by approximately 800 iterations and then remains nearly constant until around 1100 iterations.
At the same time, $\sigma_{\max}$ remains bounded and exhibits only small oscillations. The minimum control magnitude stays close to zero during most of the optimization process. These behaviors indicate that the optimization algorithm can maintain properly scaled search directions. In contrast, the formulation without thrust regularization exhibits strong oscillatory behavior throughout the optimization. Although the objective function value occasionally decreases, it does not converge to a stable value even after many iterations. The maximum singular value of the Jacobian matrix also shows large fluctuations. These behaviors indicate that the search direction is frequently damped and suggest poor local conditioning of the NLP problem. Furthermore, with thrust regularization, the minimum thrust norm value represented as $\log_{10} (\min_{i,j} \Vert{}\Delta v_{i,j} \Vert{})$ shows a trend where, as iterations progress, it smoothly and stably goes down to effectively zero. This indicates that portions of the trajectory stably settle into the optimal coasting condition. Conversely, without thrust regularization, the minimum thrust norm tends to repeatedly approach zero and then spike again. These results clearly demonstrate that the regularization approach successfully mitigates the impact of singularities associated with bang-bang control, thereby enabling robust convergence even when the optimal solution contains zero-thrust arcs.

These comparison results are consistent with the intended effect of thrust regularization. By removing the singular behavior associated with the thrust-norm representation, the regularized formulation prevents the gradient conditioning from deteriorating when the optimal solution contains near-zero-thrust arcs. As a result, the optimization algorithm can maintain stable step acceptance and convergence behavior. In contrast, the unregularized formulation remains sensitive to near-zero thrust magnitudes, which causes large variations in the conditioning of the Jacobian matrix.

\section{Numerical Results} \label{sec:example}

\subsection{Problem Setting: Mercury transfer mission inspired by the BepiColombo mission}

This study demonstrates the proposed framework using the transfer trajectory of a Mercury exploration mission inspired by BepiColombo, which employs low-thrust propulsion and nine flybys. BepiColombo is a joint ESA--JAXA mission designed to comprehensively explore Mercury and its environment using two scientific orbiters, MPO and Mio \cite{arianespace2018va245,Benkhoff2021-bc}. The BepiColombo-inspired mission considered here includes one Earth flyby, two Venus flybys, and six Mercury flybys, together with solar electric propulsion \cite{Benkhoff2021-bc,arianespace2018va245}. Because the spacecraft must reduce its relative velocity to an extremely low level to enable Mercury capture, this mission requires advanced trajectory optimization techniques involving both multiple flybys and low-thrust propulsion. Therefore, it provides an appropriate benchmark problem for demonstrating the proposed framework.

The problem setting in this study partially reflects information from actual spacecraft and launch systems. This study models the solar electric propulsion (SEP) system based on the QinetiQ T6 thruster and computes its thrust and specific impulse by interpolating publicly available performance data \cite{Hutchins2015-sb}. This study assumes Falcon Heavy Recovery as the launch vehicle~\cite{NASA_LVPerf}. Although the actual BepiColombo mission used an Ariane 5 launch vehicle, the $C_3$--payload mass capability of Falcon Heavy Recovery agrees closely with that of Ariane 5 in the energy range required for this mission. Therefore, Falcon Heavy Recovery provides a suitable substitute launch vehicle for this benchmark problem. Table~\ref{tab:model_coefficients} summarizes the coefficient sets $(\eta_{\mathrm{LV},1}, \dots,\eta_{\mathrm{LV},6})$,
$(\eta_{T,1}, \dots,\eta_{T,5})$, and $(\eta_{I_{\mathrm{sp},1}}, \dots,\eta_{I_{\mathrm{sp},5}})$ with correspond to
Eqs.~\eqref{equ:mass delivered}, \eqref{equ:thrust}, and \eqref{equ:isp}, respectively. On the other hand, detailed spacecraft engineering constraints, such as thermal conditions, navigation requirements, communication constraints, attitude constraints, and operational limitations during flybys, are not explicitly considered.

Table~\ref{tab:input_formal} summarizes the input parameters used for the trajectory search and local optimization in this example. In the local optimization, the maximum number of iterations is set to 4000 for all candidates. For each candidate, the maximum allowable total time of flight, $\tau_{\mathrm{total,max}}$, is set to 1.5 times the time of flight of the corresponding initial solution. The broad-search computations are performed on a Windows PC equipped with 32 GB of RAM and a 13th Gen Intel Core i7-13620H processor, using six CPU cores. The local-optimization computations are performed on a Mac Studio equipped with 48 GB of RAM and an Apple M4 Max processor, using three CPU cores.

\begin{table}[!htbp]
\centering
\caption{Input parameters for the trajectory search}
\label{tab:input_formal}
\begin{tabular}{lll}
\hline
\textbf{Input} & \textbf{Value} & \textbf{Comment} \\
\hline

\multicolumn{3}{l}{\textbf{For \textit{Star}-based broad search}} \\

${E}_{1-2}$ & Earth & Departure body \\
${E}_{3-4}$ & Venus & Inner-planet flybys \\
${E}_{5-11}$ & Mercury & Final target sequence \\

$t_{1}$ & (2015, 2020) & Launch window \\

$\tau_{1,11}$ & (1461, 2735) days & Total time of flight \\
$\tau_{1,2}$ & (0, 549) days & Earth--Earth transfer \\
$\tau_{i,i+1}$ & (0, 484) days & Adjacent legs \\
$\tau_{i,i+2}$ & (0, 836) days & Two-leg transfers \\

$\delta t_{i}$ & [2,2,2,2,2,2,3,4,3,3,2] days & Time mesh per encounter \\

$a_{\min,i}$ & 200 km & Minimum periapsis altitude \\

$\Delta v_{\mathrm{enc},i}$ & 2.1 km/s & Encounter $\Delta v$ allowance \\
$\Delta v_{\mathrm{total}}$ & 35 km/s & Total $\Delta v$ constraint \\

$v_{\infty,\min}$ & 0.3 km/s & Flyby lower bound \\
$v_{\infty,\max}$ & [4,9,9,9,9,9,5,4,3,0.6] km/s & Per-encounter upper bounds \\

$\Delta v_{\infty}$ & 0.05 km/s & Resonant sampling step \\

$N_{\mathrm{rev},1-2}$ & $\{0,\pm1\}$ & Low-rev transfers \\
$N_{\mathrm{rev},3-10}$ & $\{0,\pm1,\pm2,\pm3,\pm4,\pm5\}$ & Multi-rev transfers \\

$\kappa$ & $[0,\pi]$ & Crank angle candidates \\

$\Delta v_{\mathrm{leg, max}}$ & 1.0 km/s & Maximum allowable continuous thrust $\Delta v$ per leg \\

$\zeta$ & 1.5 & Threshold ratio for triggering convex pruning \\

$N_j$ & 31 & Number of control nodes within the linearized MEE grid \\

$a_{\mathrm{acc},\mathrm{max}}$ & $6.525 \times 10^{-8}$ $\mathrm{km/s}^2$ & Maximum thrust acceleration limit ($D \cdot N_{\mathrm{EP}} T_{\mathrm{thrust},\mathrm{max}} / m_{\mathrm{ref}}$ with $m_{\mathrm{ref}} = 4000$ kg) \\
\hline

\hline

\multicolumn{3}{l}{\textbf{For Local Optimization}} \\

$\sigma_\mathrm{LV}$ & 0.1 & Launch performance margin \\
${C_3}_\mathrm{max}$ & 12.1 $\mathrm{km}^2/\mathrm{s}^2$ & Max launch $C_3$  \\
$D$ & 0.9 & Duty cycle \\
$m_{\mathrm{dry}}$ & 3305.5 kg & Dry mass \\
$m_\mathrm{max}$ & 4043.2 kg & Maximum initial mass \\
$P_{\mathrm{BOL}}$ & 14 kW & Beginning-of-life power \\
$\delta _\mathrm{power}$ & 0.15 & Power margin \\
$\tau_\mathrm{deg}$ & 0.01 & Degradation factor \\
$N_{\mathrm{EP}}$ & 2 & Number of thrusters \\
$T_{\mathrm{thrust},\mathrm{max}}$ & $145$ mN & Maximum thrust performance per engine unit \\
$I_{\mathrm{sp},\mathrm{max}}$ & 4300 s & Maximum engine specific impulse capability\\
$N_j$ & $21 (|N_{\mathrm{rev},i}| + 1)$ & Number of SFT nodes for leg $i$\\
\hline
\end{tabular}
\end{table}

\begin{figure}
    \centering
    \includegraphics[width=1.0\linewidth]{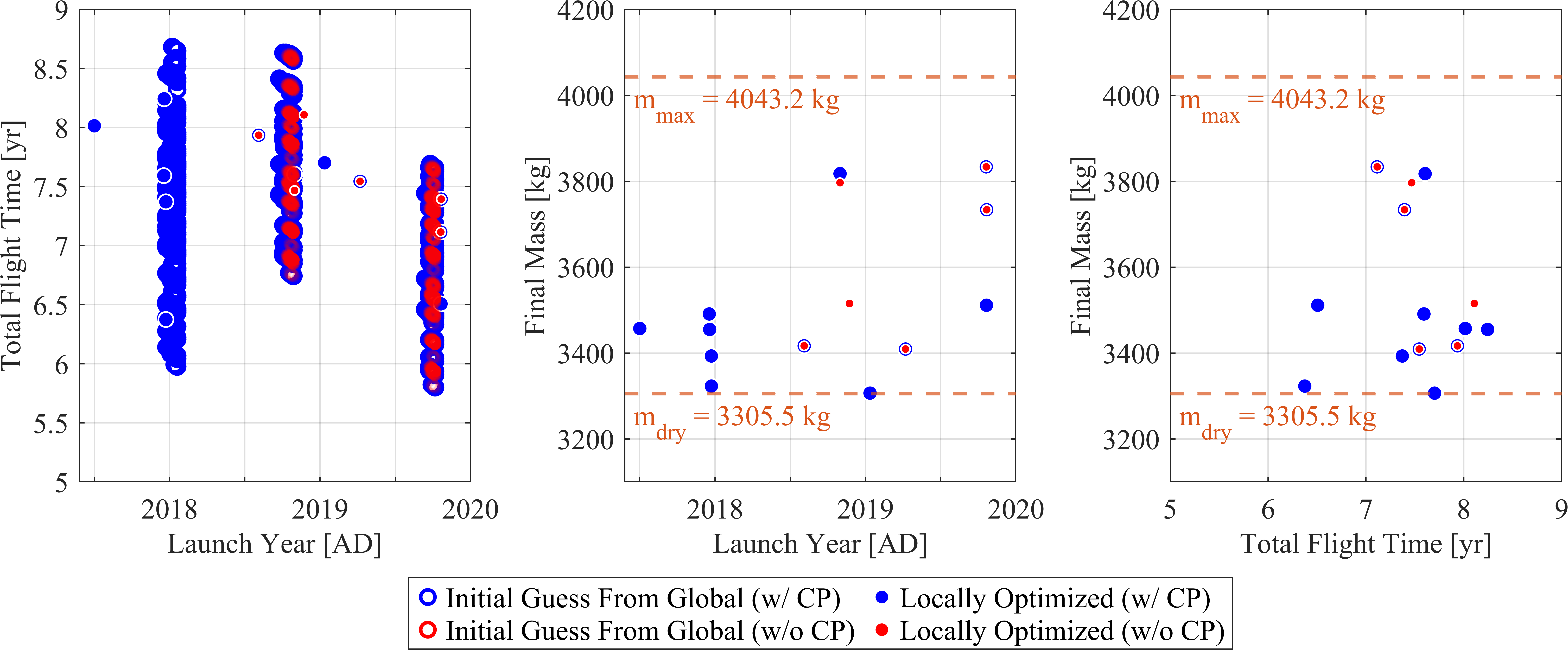}
    \caption{
    (a) Launch-year--flight-time solution space; final-mass trade-offs with (b) launch year and (c) total flight time before and after optimization, with and without convex pruning.}
    \label{fig:Trade-Off}
\end{figure}

\subsection{Broad Search}
 To assess the effect of the CP, this study performs the broad search both with and without CP and compares the resulting solution sets after local optimization. Without CP, hereafter denoted as w/o CP, the broad search requires 1578.563 sec, of which 584.579 sec is spent in the \textit{Leg} Module, and generates 4,664 feasible candidate solutions.
 The framework clusters the candidates obtained from the broad search into one-week launch-date intervals and ranks them according to $\Delta v_\mathrm{total}$. The differences in computation time and the number of obtained candidate solutions are examined first. These candidates are classified into six clusters corresponding to one-week launch windows. The fractions of the total runtime spent in the \textit{Body}, \textit{Leg}, \textit{Flyby}, and \textit{Combo} modules were 0.20\%, 37.92\%, 55.00\%, and 6.88\%, respectively. With CP, hereafter denoted as w/ CP, the total runtime increased to 4915.74 sec, of which 3113.623 sec is spent in the \textit{Leg} Module. However, the number of feasible candidate solutions increases to 36335, and the number of launch-window clusters increases to 15. In this case, the fractions of the total runtime spent in the \textit{Body}, \textit{Leg}, \textit{Flyby}, and \textit{Combo} modules are 0.06\%, 63.93\%, 21.42\%, and 14.59\%, respectively. These results show that, although CP introduces an additional computational cost, it substantially expands the candidate-solution search space within a practical computational range.

The broad search results in Fig.~\ref{fig:Trade-Off} (a) provide a visual comparison of the candidate solutions obtained in the w/ CP and w/o CP cases in terms of launch year and total flight time. The distribution of the candidate solutions shows that the w/ CP candidates cover the design space obtained in the w/o CP case while increasing both the number of candidates and the number of clusters. In particular, the proposed method obtains families of initial solutions in the departure periods around mid-2017 and early 2018, which are not found in the w/o CP case. In addition, it identifies a larger number of candidate families for departure periods around late 2018, which corresponds to the actual BepiColombo launch period, and from mid- to late 2019.
These results indicate that CP recovers candidates that would otherwise be discarded and substantially expands the feasible design space.

\subsection{Local Optimization}

The differences in computation time for local optimization are examined using the broad search candidate solutions obtained with and without CP. The local optimization is performed on the same computer using three CPU cores. For converged cases, the computation time is generally about 300--800 sec per case. Non-converged cases terminate after reaching the maximum number of iterations and require about 800--900 s. For the candidates obtained without CP, all six selected clusters yield converged solutions, with an average computation time of $3.024~\mathrm{h}$ per cluster. For the candidates obtained with CP, 12 of the 15 clusters yield converged solutions, with an average computation time of $1.349~\mathrm{h}$ per cluster. The average time required for convergence or for testing all candidates in each cluster is 55.4\% lower than that for the candidates obtained without CP. The three non-converged clusters contain only 8, 12, and 12 representative candidates, respectively, so all available candidates are exhausted before a converged solution can be obtained. 

Table~\ref{tab:rank_summary} evaluates the effect of the initial control histories generated by CP on this reduction in computation time for local optimization. Of the 12 converged clusters, eight encounter at least one candidate with a CP-generated initial control history before convergence, whereas the remaining four converge before any such candidate is evaluated. Additional convergence tests therefore examine these eight clusters. For each cluster, the table compares the computation time required for convergence between cases with and without the CP-generated initial control histories. The CP-generated initial control histories reduce the mean computation time from 2.68 h to 1.48 h, corresponding to a reduction of approximately 44.8\%.  They also reduce the total computation time from 21.45 h to 11.85 h. These results suggest that the control histories generated by CP improve the convergence behavior of low-thrust trajectory optimization by providing favorable initial search directions and thereby shorten the computation time required to obtain converged solutions.
 However, the improvement is not observed for all trajectory candidates. In some clusters, the cases without initial control histories converge more quickly. This behavior suggests that, in some cases, an almost ballistic control profile may be closer to the local optimal solution.

 \begin{table}[hbt!]
    \centering
    \caption{Summary of single-core computation times for successfully converged clusters, comparing cases with and without the CP-based initial control guess.}
    \label{tab:rank_summary}
    \begin{tabular}{lcc}
        \hline\hline
        & w/ Initial Control Guess & w/o Initial Control Guess \\
        \hline
        Mean Time [h] & 1.48 & 2.68 \\
        Total Computation Time [h] & 11.85 & 21.45 \\
        \hline\hline
    \end{tabular}
\end{table}

Fig.~\ref{fig:Trade-Off} (b) and (c) show that the use of CP expands the design space covered by the final locally optimized trajectories. A particularly notable result appears in the departure window around early 2018. This region is not identified by the broad search without CP, but it produces the optimized trajectory with the shortest total flight time among the obtained solutions. As a result, the proposed framework expands the departure-epoch range of the optimized solutions by more than one year. The relationship between final mass and total flight time also shows an important trend. Although the trade-off is not strictly monotonic, longer transfer times generally tend to produce larger final masses. This tendency is consistent with the characteristics of low-thrust propulsion, because a longer flight time provides more opportunities to use low-thrust acceleration efficiently and reduce propellant consumption. These local-optimization results show that the proposed pruning method captures practically important trajectory families that could be missed by conventional broad search without low-thrust effects and expands the search space of optimized solutions.

\subsection{Representative Optimized Trajectories} \label{sec:representative}

\begin{sidewaystable}[htbp]
\centering
\caption{Summary of optimal solutions for the BepiColombo-inspired scenario mission with and without convex pruning.}
\label{tab:optimal_design_summary}
\small
\setlength{\tabcolsep}{3pt}
\renewcommand{\arraystretch}{1.5}
\begin{tabular}{llccccccccccccc}
\hline
Group & ID & $\Delta V_{\mathrm{GA},1}$ & $\Delta V_{\mathrm{GA},2}$ & $\Delta V_{\mathrm{GA},3}$ & $\Delta V_{\mathrm{GA},4}$ & $\Delta V_{\mathrm{GA},5}$ & $\Delta V_{\mathrm{GA},6}$ & $\Delta V_{\mathrm{GA},7}$ & $\Delta V_{\mathrm{GA},8}$ & $\Delta V_{\mathrm{GA},9}$ & Final Mass [kg] & Launch Date & Arrival Date & TOF [yr] \\
\hline
Optimal Trajectory (w/ CP) & 1 & -3.00 & -1.71 & -6.85 & -1.36 & -1.38 & -0.80 & -2.86 & -0.10 & -0.19 & 3393.1 & 12/23/2017 & 05/07/2025 & 7.37 \\ 
 & 2 & -2.98 & -1.73 & -6.96 & -1.27 & -2.08 & -2.49 & -0.26 & -0.02 & -0.01 & 3306.6 & 01/12/2019 & 09/25/2026 & 7.70 \\ 
 & 3 & -3.00 & -1.71 & -6.85 & -1.36 & 0.02 & -0.03 & -1.94 & -2.74 & -0.04 & 3457.1 & 07/02/2017 & 07/07/2025 & 8.01 \\ 
 & 4 & -2.46 & -1.74 & -7.14 & -1.20 & -1.37 & -3.53 & -0.22 & -0.03 & -0.02 & 3511.3 & 10/21/2019 & 04/24/2026 & 6.51 \\ 
 & 5 & -3.00 & -1.71 & -6.85 & -1.36 & 0.00 & -2.05 & -2.61 & 0.07 & -0.17 & 3416.8 & 08/05/2018 & 07/12/2026 & 7.93 \\ 
 & 6 & -2.98 & -1.72 & -6.85 & -1.36 & -1.37 & -0.78 & -2.93 & 0.19 & -0.19 & 3817.6 & 10/31/2018 & 06/08/2026 & 7.60 \\ 
 & 7 & -2.46 & -1.74 & -7.14 & -1.22 & -2.25 & 0.70 & -3.47 & 0.13 & -0.16 & 3733.6 & 10/21/2019 & 03/14/2027 & 7.39 \\ 
 & 8 & -2.48 & -1.71 & -6.85 & -1.36 & -1.38 & -0.00 & -0.80 & -3.05 & -0.01 & 3409.3 & 04/08/2019 & 10/23/2026 & 7.54 \\ 
 & 9 & -3.00 & -1.71 & -6.85 & -1.36 & -1.37 & -0.78 & -2.93 & 0.19 & -0.19 & 3491.0 & 12/18/2017 & 07/22/2025 & 7.59 \\ 
 & 10 & -2.98 & -1.73 & -6.96 & -1.27 & -2.08 & -2.49 & -0.30 & 0.01 & -0.02 & 3323.0 & 12/23/2017 & 05/08/2024 & 6.37 \\ 
 & 11 & -3.00 & -1.71 & -6.85 & -1.36 & -1.35 & -0.02 & -3.27 & 0.02 & -0.15 & 3454.9 & 12/19/2017 & 03/17/2026 & 8.24 \\ 
 & 12 & -2.46 & -1.74 & -7.14 & -1.22 & -1.41 & -0.01 & -3.88 & 0.16 & -0.21 & 3833.4 & 10/20/2019 & 12/01/2026 & 7.12 \\ 
\hline
Optimal Trajectory (w/o CP) & 1 & -2.46 & -1.74 & -7.14 & -1.17 & -0.04 & -1.97 & -2.78 & 0.17 & -0.19 & 3733.6 & 10/21/2019 & 03/14/2027 & 7.39 \\ 
 & 2 & -2.46 & -1.74 & -7.14 & -0.06 & -2.54 & -0.76 & -2.77 & 0.17 & -0.16 & 3833.4 & 10/20/2019 & 12/01/2026 & 7.12 \\ 
 & 3 & -2.48 & -1.71 & -6.85 & -1.36 & -1.38 & -0.00 & -0.80 & -3.05 & -0.01 & 3409.3 & 04/08/2019 & 10/23/2026 & 7.54 \\ 
 & 4 & -2.98 & -1.72 & -6.92 & -1.32 & -1.36 & -3.52 & -0.22 & -0.03 & -0.02 & 3796.5 & 10/31/2018 & 04/19/2026 & 7.47 \\ 
 & 5 & -3.00 & -1.71 & -6.77 & -1.41 & -0.03 & -2.04 & -2.61 & 0.07 & -0.17 & 3416.8 & 08/05/2018 & 07/12/2026 & 7.93 \\ 
 & 6 & -2.98 & -1.72 & -6.92 & -1.32 & -1.36 & -0.77 & -2.72 & 0.10 & -0.21 & 3515.4 & 11/23/2018 & 01/01/2027 & 8.11 \\ 
\hline
\end{tabular}

\end{sidewaystable}

Table~\ref{tab:optimal_design_summary} provides more detailed information on the optimized solutions obtained with and without CP.
The table identifies each optimized solution by its ID and lists its performance metrics, including the heliocentric speed gain or loss induced by each gravity assist ($\Delta v_\mathrm{GA,i}$), where $i$ denotes the gravity-assist index, the final mass, the departure and arrival epochs, and the total time of flight.
In this comparison, all optimized solutions reach the prescribed launch-energy upper bound, $C_3 = 12.1~\mathrm{km^2/s^2}$, and therefore also reach the allowable upper bound on the initial mass. The table shows that both solution sets contain trajectories with relatively large final mass and short flight time. For example, the w/ CP solution with ID 12 and the w/o CP solution with ID 2 correspond to the almost identical trajectory and achieve the largest final mass among the listed solutions. This result indicates that trajectories requiring relatively small low-thrust corrections can already be identified by the broad search without CP. This tendency is consistent with the fact that, when a solution requires less low-thrust acceleration, each leg tends to contain longer coast arcs and can therefore be approximated sufficiently well by a search that does not explicitly account for low-thrust effects. On the other hand, the w/ CP solution with ID 10 achieves the shortest total flight time among all listed solutions. Its total flight time is 6.37 yr, which is approximately 0.75 yr shorter than the shortest solution obtained without CP. This result suggests that, when the feasible time of flight around the initial trajectory is short, the trajectory requires a larger contribution from low-thrust acceleration to satisfy the prescribed transfer time. Such solutions may not appear in the local-optimization stage unless the broad search identifies suitable local optima basins from the combinatorial solution space while accounting for low-thrust corrections through CP.

\begin{figure}[t!]
    \centering
    
    \begin{minipage}{0.48\linewidth} 
        \centering
        \includegraphics[width=\linewidth]{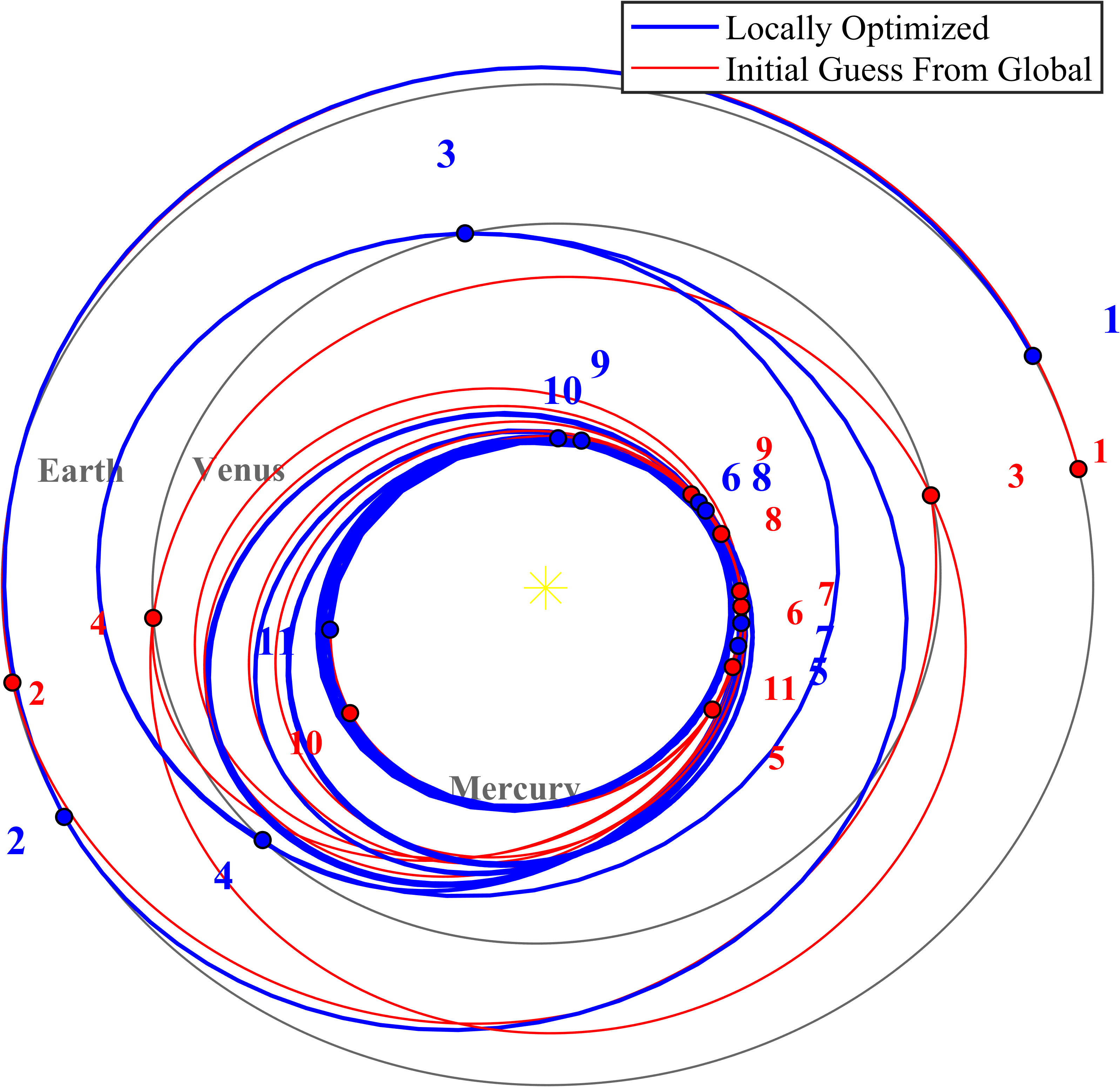}
    \end{minipage}
    \hfill 
    \begin{minipage}{0.48\linewidth}
        \centering
        \includegraphics[width=\linewidth]{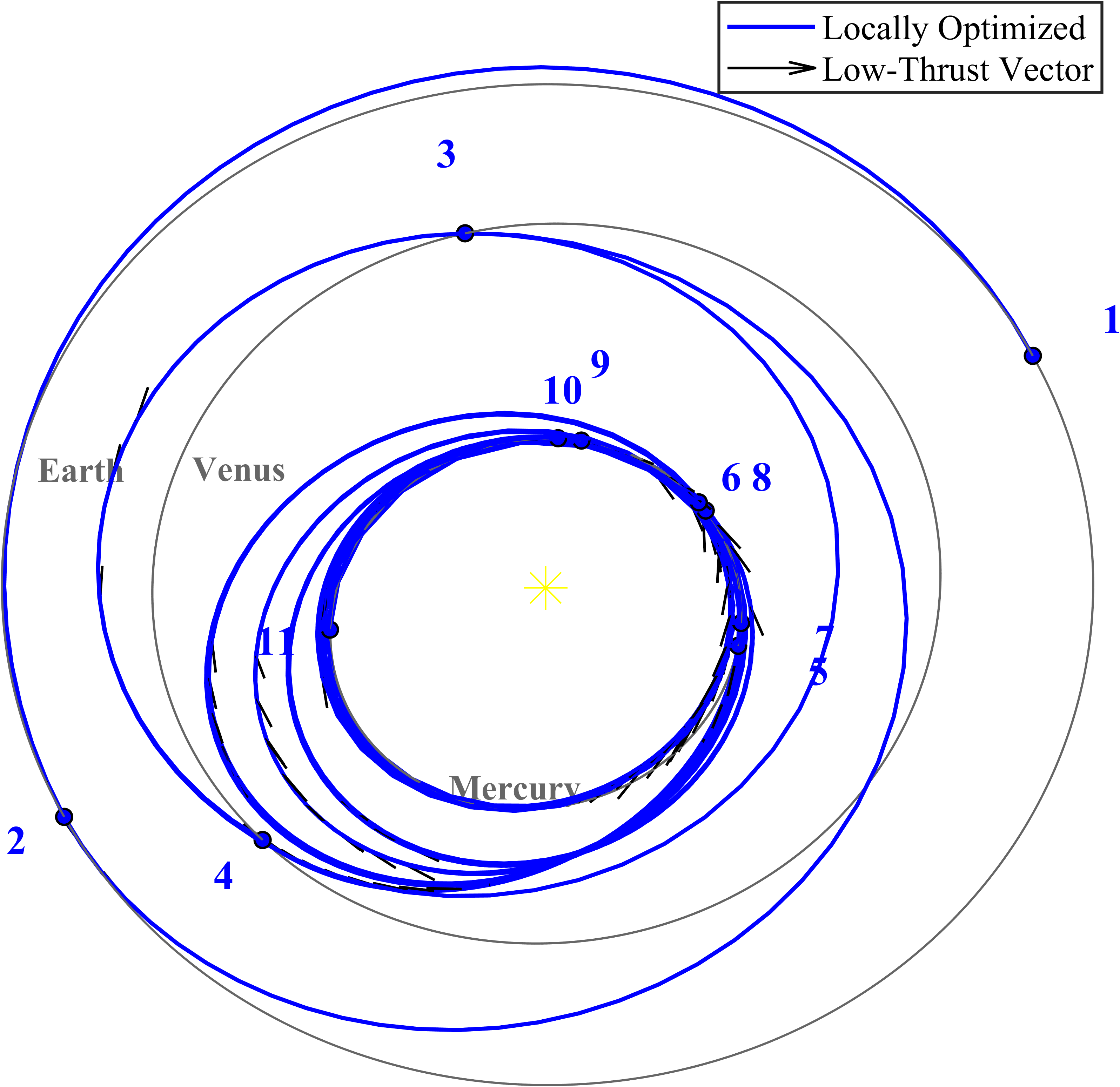}
    \end{minipage}
    
    \vspace{0.5cm} 
    
    \includegraphics[width=1.0\linewidth]{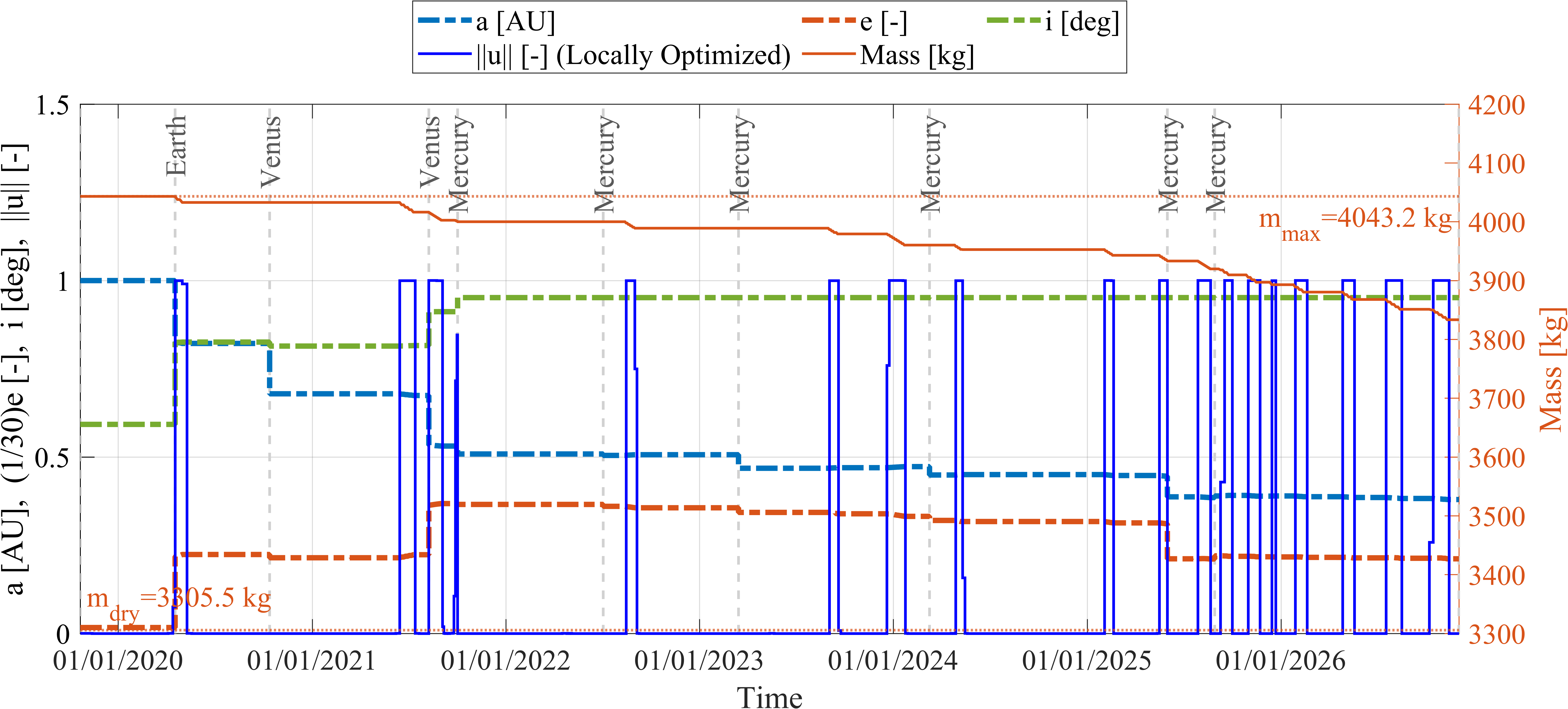}

    \caption{Optimized trajectory, thrust vectors, and time-history profiles for the maximum-final-mass representative solution, corresponding to the w/ CP solution with ID 12.}
    \label{fig:bepi_plots_maxmf}
\end{figure}

\begin{table}[tb]
\centering
\caption{Mission itinerary and event details for the maximum-final-mass solution obtained with CP, ID 12: where the locally optimized solutions are shown in bold.}
\label{tab:bepi_events_maxmf}   
\small
\setlength{\tabcolsep}{4pt}
\renewcommand{\arraystretch}{1.7}
\begin{tabular}{llccccc}
\hline
Event & Date & $ v_{\infty}$ (km/s) & Altitude (km) & Mass (kg) & $\Delta v_{\mathrm{Mismatch}}$ (km/s) & $\Delta v_{\mathrm{Leg}}$ (km/s) \\
\hline
1  Earth & 10/07/2019 $|$ \textbf{10/21/2019} & 3.550 $|$ \textbf{3.475} & -- $|$ \textbf{--} & -- $|$ \textbf{4043} & -- $|$ \textbf{--} & 0.000 $|$ \textbf{0.004} \\ 
2  Earth & 03/31/2020 $|$ \textbf{04/17/2020} & 3.550 $|$ \textbf{3.451} & {37330} $|$ \textbf{4568} & -- $|$ \textbf{4043} & 2.066 $|$ \textbf{0.000} & 0.000 $|$ \textbf{0.099} \\ 
3  Venus & 08/18/2020 $|$ \textbf{10/12/2020} & 6.666 $|$ \textbf{7.408} & 200 $|$ \textbf{12041} & -- $|$ \textbf{4033} & 1.581 $|$ \textbf{0.000} & 0.000 $|$ \textbf{0.164} \\ 
4  Venus & 07/14/2021 $|$ \textbf{08/08/2021} & 8.247 $|$ \textbf{7.652} & {533} $|$ \textbf{200} & -- $|$ \textbf{4016} & 1.290 $|$ \textbf{0.000} & 0.000 $|$ \textbf{0.163} \\ 
5  Mercury & 09/26/2021 $|$ \textbf{10/01/2021} & 6.350 $|$ \textbf{6.556} & {706} $|$ \textbf{200} & -- $|$ \textbf{4000} & 0.969 $|$ \textbf{0.000} & 0.000 $|$ \textbf{0.000} \\ 
6  Mercury & 06/25/2022 $|$ \textbf{07/03/2022} & 5.124 $|$ \textbf{6.073} & 200 $|$ \textbf{36309} & -- $|$ \textbf{4000} & 1.129 $|$ \textbf{0.000} & 0.000 $|$ \textbf{0.110} \\ 
7  Mercury & 03/17/2023 $|$ \textbf{03/15/2023} & 3.921 $|$ \textbf{5.705} & {3773} $|$ \textbf{200} & -- $|$ \textbf{3989} & 0.811 $|$ \textbf{0.000} & 0.000 $|$ \textbf{0.289} \\ 
8  Mercury & 03/07/2024 $|$ \textbf{03/09/2024} & 2.920 $|$ \textbf{3.616} & {4262} $|$ \textbf{200} & -- $|$ \textbf{3960} & 2.015 $|$ \textbf{0.000} & 0.000 $|$ \textbf{0.275} \\ 
9  Mercury & 05/24/2025 $|$ \textbf{05/31/2025} & 0.809 $|$ \textbf{2.861} & {10069 } $|$ \textbf{200} & -- $|$ \textbf{3933} & 0.204 $|$ \textbf{0.000} & 0.000 $|$ \textbf{0.137} \\ 
10  Mercury & 09/25/2025 $|$ \textbf{08/28/2025} & 0.437 $|$ \textbf{2.598} & {14651} $|$ \textbf{66542} & -- $|$ \textbf{3920} & 0.204 $|$ \textbf{0.000} & 0.000 $|$ \textbf{0.893} \\ 
11  Mercury & 04/28/2026 $|$ \textbf{12/02/2026} & 0.432 $|$ \textbf{0.500} & -- $|$ \textbf{--} & -- $|$ \textbf{3833} & -- $|$ \textbf{--} & -- $|$ \textbf{--} \\ 
\hline
\end{tabular}

\end{table}

\begin{figure}[t!]
    \centering
    
    \begin{minipage}{0.48\linewidth} 
        \centering
        \includegraphics[width=\linewidth]{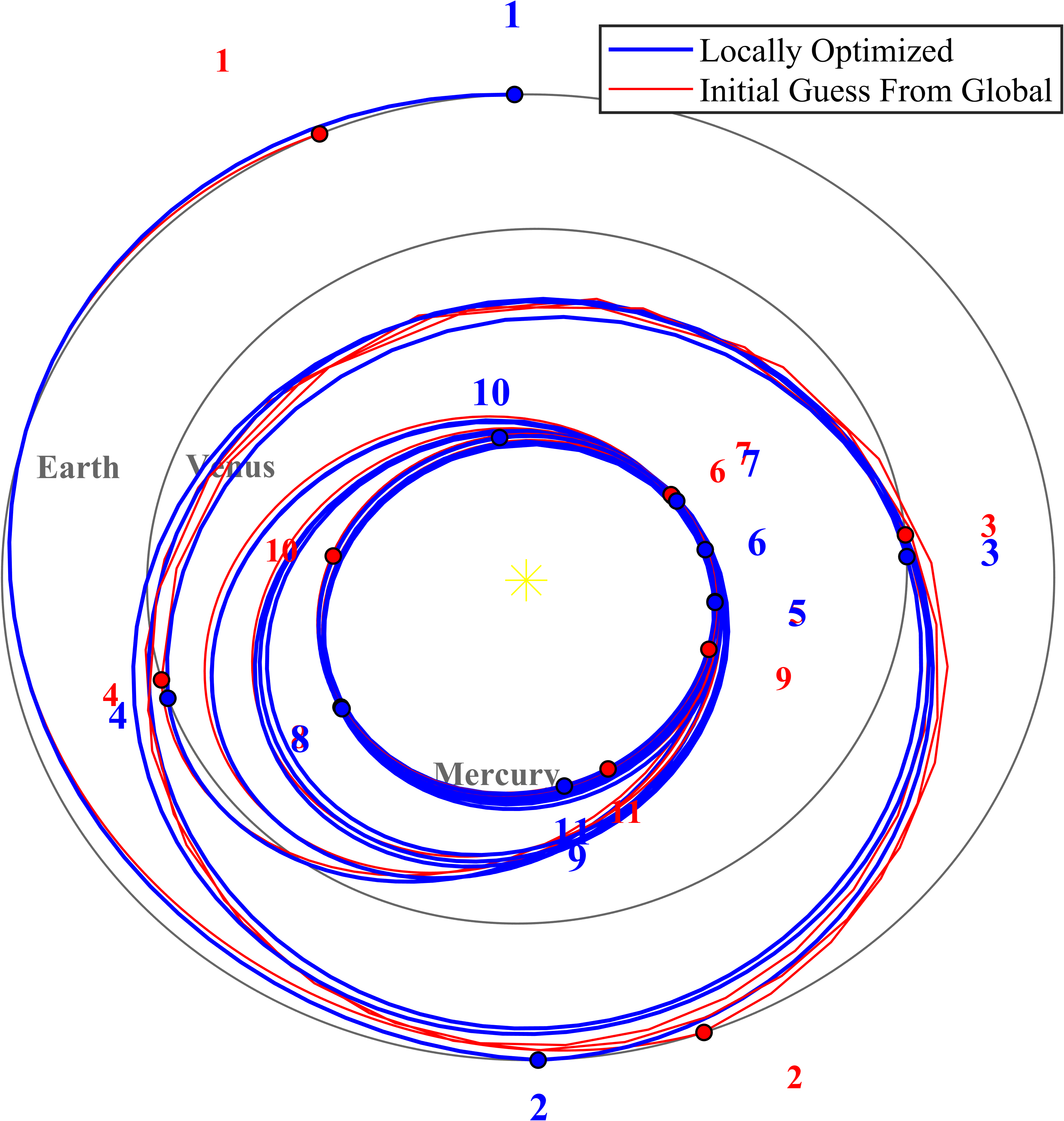}
    \end{minipage}
    \hfill 
    \begin{minipage}{0.48\linewidth}
        \centering
        \includegraphics[width=\linewidth]{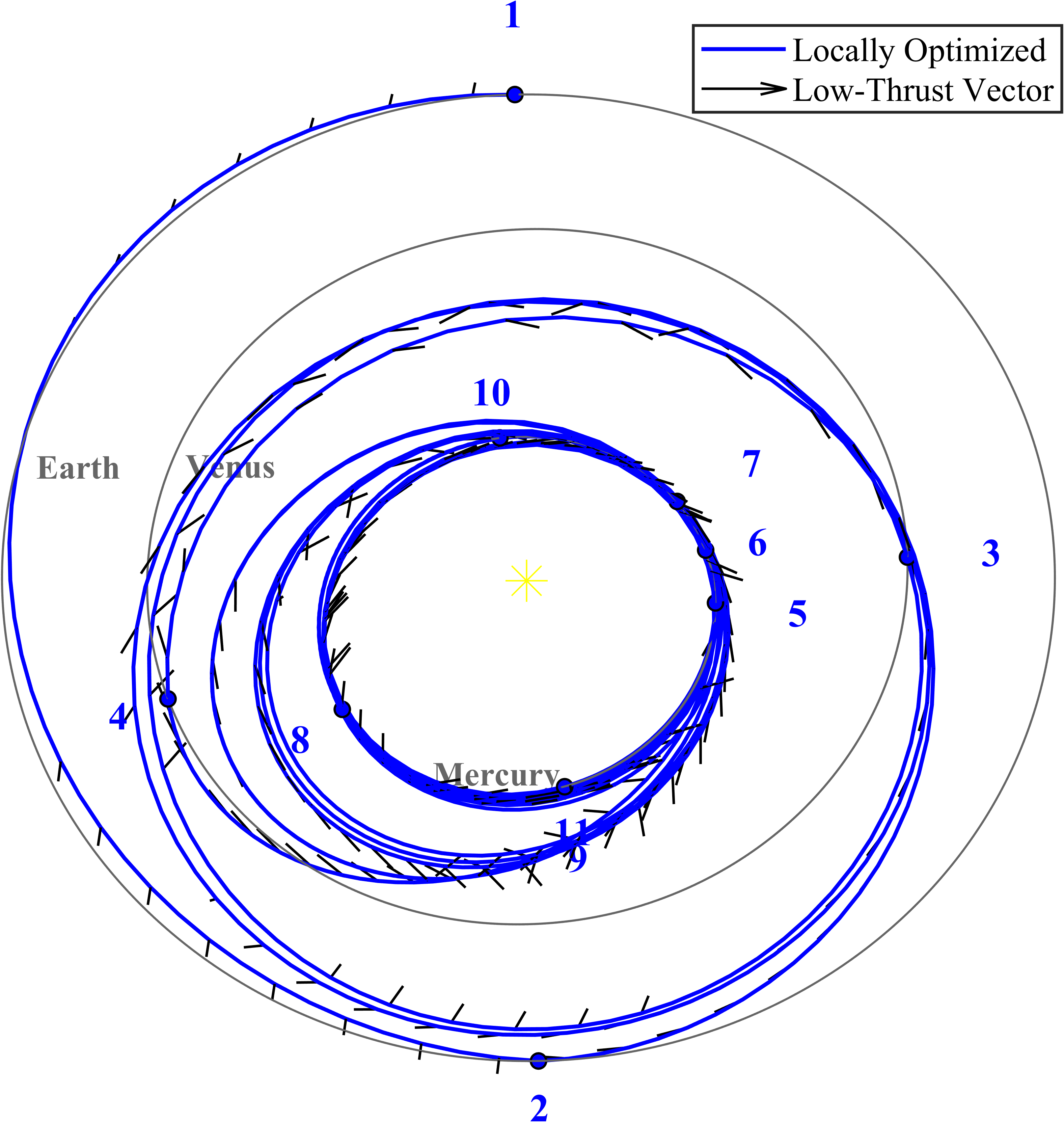}
    \end{minipage}
    
    \vspace{0.5cm} 
    
    \includegraphics[width=1.0\linewidth]{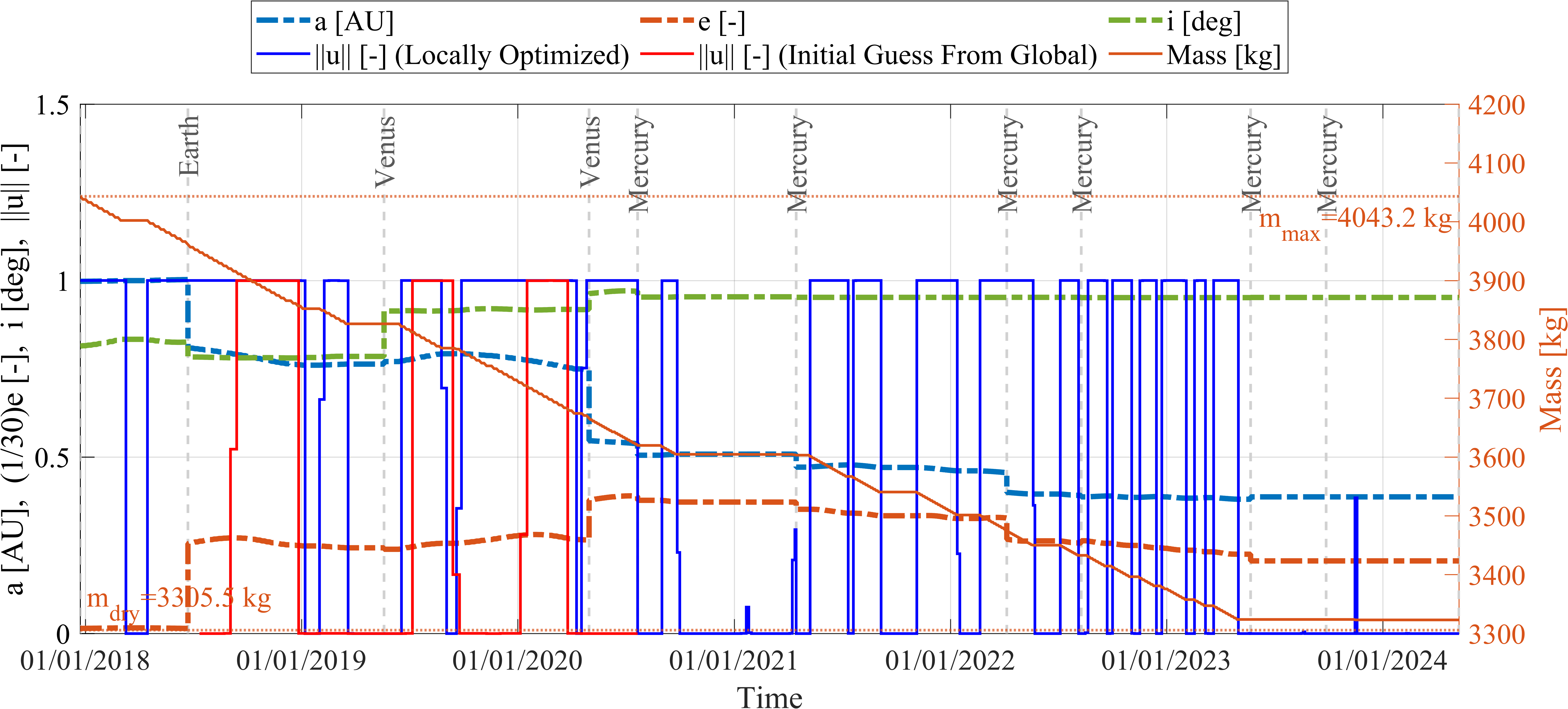}

        \caption{Optimized trajectory, thrust vectors, and time-history profiles for the minimum-time representative solution, corresponding to the w/ CP solution with ID 10.}
    \label{fig:bepi_plots_mintof}
\end{figure}

Tables~\ref{tab:bepi_events_maxmf} and \ref{tab:bepi_events_mintof} summarize the detailed mission itineraries of the maximum-final-mass representative solution, corresponding to the w/ CP solution with ID~12, and the minimum-time-of-flight representative solution, corresponding to the w/ CP solution with ID~10, respectively. For every encounter, the tables report the encounter epoch, hyperbolic-excess velocity ($v_\infty$), flyby altitude, spacecraft mass, the flyby-mismatch $\Delta v$ ($\Delta v_\mathrm{mismatch}$), and low-thrust expenditure over the subsequent leg ($\Delta v_\mathrm{Leg}$) in the initial and optimized solutions.
Figures~\ref{fig:bepi_plots_maxmf} and \ref{fig:bepi_plots_mintof} present the corresponding trajectory and time-history results. The top left panel compares the initial and optimized trajectories, the top right panel shows the optimized trajectory together with the low-thrust vectors, and bottom panel provides the histories of the orbital elements $a$, $e$, and $i$, low-thrust control, and spacecraft mass. The numbers in the top panel indicate the encounter sequence.

This analysis first examines the maximum-final-mass solution. As shown in Table~\ref{tab:optimal_design_summary}, the spacecraft reaches Mercury with a final mass of $3833.4~\mathrm{kg}$, which corresponds to approximately 28\% of the propellant consumed from the original mission. The time of flight is $7.12~\mathrm{yr}$, which is approximately $0.97~\mathrm{yr}$ shorter than the reference transfer duration of BepiColombo. The thrust profile in Fig.~\ref{fig:bepi_plots_maxmf} exhibits a bang-bang-like control structure consistent with the optimality condition from Pontryagin's maximum principle, while many legs also contain coast arcs. This behavior is consistent with the tendency that such a solution can be identified as feasible even without a CP-generated initial control history shown in Table~\ref{tab:bepi_events_maxmf}, because low-thrust acceleration plays a relatively limited role.
In addition, the top left panel of Fig.~\ref{fig:bepi_plots_maxmf} shows a relatively large deviation between the initial and optimized trajectories for ID 12, which does not have a CP-generated initial control history. By contrast, the top left panel of Fig.~\ref{fig:bepi_plots_mintof} shows a much smaller deviation for ID 10, whose initial solution includes a control history generated by CP. This contrast suggests that the CP provides an initial trajectory and a control guess that is closer to the locally optimized low-thrust solution.
The histories of the orbital elements $a$, $e$, and $i$ show large variations at the flyby epochs. This behavior indicates that gravity assists produce most of the trajectory changes, whereas low thrust mainly provides fine adjustments. Therefore, compared with the reference mission, this trajectory efficiently uses gravity assists, preserves most of the propellant mass, and achieves an earlier arrival at Mercury.

\begin{table}[tb]
\centering
\caption{Mission itinerary and event details for the minimum-time solution obtained with CP, ID 10: where the locally optimized solutions are shown in bold.}
\label{tab:bepi_events_mintof}   
\small
\setlength{\tabcolsep}{4pt}
\renewcommand{\arraystretch}{1.7}
\begin{tabular}{llccccc}
\hline
Event & Date & $ v_{\infty}$ (km/s) & Altitude (km) & Mass (kg) & $\Delta v_{\mathrm{Mismatch}}$ (km/s) & $\Delta v_{\mathrm{Leg}}$ (km/s) \\
\hline
1  Earth & 01/13/2018 $|$ \textbf{12/23/2017} & 3.950 $|$ \textbf{3.475} & -- $|$ \textbf{--} & -- $|$ \textbf{4043} & -- $|$ \textbf{--} & 0.000 $|$ \textbf{0.820} \\ 
2  Earth & 07/12/2018 $|$ \textbf{06/23/2018} & 4.950 $|$ \textbf{3.973} & {1080} $|$ \textbf{581} & -- $|$ \textbf{3961} & 0.017 $|$ \textbf{0.000} & 0.621 $|$ \textbf{1.387} \\ 
3  Venus & 05/22/2019 $|$ \textbf{05/20/2019} & 9.009 $|$ \textbf{8.356} & {51404} $|$ \textbf{22080} & -- $|$ \textbf{3826} & 1.287 $|$ \textbf{0.000} & 0.800 $|$ \textbf{1.718} \\ 
4  Venus & 04/28/2020 $|$ \textbf{04/30/2020} & 10.069 $|$ \textbf{9.387} & 200 $|$ \textbf{200} & -- $|$ \textbf{3666} & 0.755 $|$ \textbf{0.000} & 0.000 $|$ \textbf{0.507} \\ 
5  Mercury & 07/21/2020 $|$ \textbf{07/21/2020} & 5.911 $|$ \textbf{6.388} & 200 $|$ \textbf{200} & -- $|$ \textbf{3620} & 2.069 $|$ \textbf{0.000} & 0.000 $|$ \textbf{0.184} \\ 
6  Mercury & 04/19/2021 $|$ \textbf{04/15/2021} & 3.444 $|$ \textbf{5.072} & {2208} $|$ \textbf{200} & -- $|$ \textbf{3603} & 1.384 $|$ \textbf{0.000} & 0.000 $|$ \textbf{1.453} \\ 
7  Mercury & 04/06/2022 $|$ \textbf{04/05/2022} & 2.064 $|$ \textbf{3.624} & 200 $|$ \textbf{200} & -- $|$ \textbf{3475} & 1.647 $|$ \textbf{0.000} & 0.000 $|$ \textbf{0.487} \\ 
8  Mercury & 08/09/2022 $|$ \textbf{08/09/2022} & 0.426 $|$ \textbf{3.009} & 200 $|$ \textbf{44395} & -- $|$ \textbf{3433} & 0.107 $|$ \textbf{0.000} & 0.000 $|$ \textbf{1.293} \\ 
9  Mercury & 03/12/2023 $|$ \textbf{05/22/2023} & 0.318 $|$ \textbf{0.509} & 200 $|$ \textbf{33093} & -- $|$ \textbf{3324} & 0.029 $|$ \textbf{0.000} & 0.000 $|$ \textbf{0.000} \\ 
10  Mercury & 07/14/2023 $|$ \textbf{09/27/2023} & 0.326 $|$ \textbf{0.749} & 200 $|$ \textbf{25182951} & -- $|$ \textbf{3324} & 0.029 $|$ \textbf{0.000} & 0.000 $|$ \textbf{0.009} \\ 
11  Mercury & 02/14/2024 $|$ \textbf{05/08/2024} & 0.330 $|$ \textbf{0.500} & -- $|$ \textbf{--} & -- $|$ \textbf{3323} & -- $|$ \textbf{--} & -- $|$ \textbf{--} \\ 
\hline
\end{tabular}

\end{table}

This analysis next examines the minimum-time-of-flight solution. As shown in Table~\ref{tab:optimal_design_summary}, this solution reaches Mercury in $6.37~\mathrm{yr}$ with a final mass of $3323.0~\mathrm{kg}$. The corresponding propellant consumption is approximately 98\% of the reference mission's propellant capacity, $737.7~\mathrm{kg}$, and the transfer duration is shortened by approximately $1.72~\mathrm{yr}$. The thrust profile of the minimum-time-of-flight solution also shows a bang-bang-like control structure. However, the active thrusting periods are much longer than those of the maximum-final-mass solution. The low-thrust distribution before 2020 indicates that the final optimized control history retains key features of the initial control history generated by CP. As shown in Table~\ref{tab:bepi_events_mintof}, the legs associated with encounter~2 at Earth,  encounter~3 at Venus, and encounter~4 at Venus require relatively large values of $\Delta v_{\mathrm{leg}}$ in both the initial and optimized solutions. This persistence suggests that the subsequent nonlinear optimization preserves the principal thrusting phases introduced by CP. In addition, the optimized hyperbolic-excess velocity $v_\infty$ increases sharply from 3.973 km/s at encounter 2 with Earth to 8.356 km/s and 9.387 km/s at encounters 3 and 4 with Venus, respectively. This rapid increase indicates a substantial change in the heliocentric transfer energy during this phase of the trajectory.
These changes occur during a phase of prolonged low-thrust operation and coincide with pronounced variations in the orbital elements $a$, $e$, and $i$ between successive legs. Because such drastic orbital shaping requires substantial energy, purely impulsive transfer models inevitably yield excessive $\Delta v$ requirements, which would normally cause these candidates to be pruned during the broad search. These observations are consistent with the intended role of CP in the broad search: CP recovers near-feasible candidate trajectories by introducing the minimum low-thrust correction required by the convexified model and provides an initial control history that remains informative during the subsequent nonlinear optimization. It also suggests that the CP can generate an initial control history whose qualitative trend is similar to that of the final optimized solution. This solution represents a fast transfer trajectory obtained by fully exploiting both low-thrust propulsion and gravity assists. The result suggests that, for finding such trajectories and obtaining stable locally optimized solutions, the proposed framework is effective because it incorporates the CP-based recovery process into the broad search stage and uses the resulting low-thrust profile in the local optimization.


\section{Discussion} \label{sec:discussion}

This study applies the proposed design framework to a BepiColombo-inspired Mercury transfer mission involving low-thrust propulsion and multiple flybys and demonstrates its performance. The results show that CP discovers trajectory candidates that cannot be captured by a lightweight broad-search model alone and expands the searchable design space. They also show that the initial trajectory guesses obtained from the broad search and the CP-generated initial control histories improve the convergence behavior of the subsequent trajectory optimization when they are combined with the thrust-regularized local-optimization model.

Although the broad-search stage in this study handles a specific constrained multivariable problem associated with low-thrust optimization, evaluating the robustness of the proposed framework under different transfer strategies and constraint sets remains an important research topic. 
The flexibility of the convex formulation suggests that the proposed method can be extended beyond the specific low-thrust optimization problem considered in this study. If the corresponding correction problem can be formulated or approximated as a convex problem, the same concept is applicable to multiple deep-space maneuvers~\cite{Vavrina2016-fy}, solar-sail trajectories~\cite{Oguri2025-hy}, boundary-velocity constraints, and thrust-direction or thrust-magnitude constraints. Applying the framework to such problems would clarify how robustly CP works under different types of trajectory corrections and mission constraints.

Reducing the computational cost of broad search with CP is another important future direction. In the BepiColombo-inspired application, the implementation with CP requires more computation time, and the CP-related process accounts for a large fraction of the total broad-search cost. In multi-flyby trajectory design, the computational burden can increase rapidly as the number of flyby sequences grows, which may limit the range of practical applications within realistic computation time. Therefore, future work should improve the computational efficiency of the framework from both algorithmic and architectural viewpoints. Possible directions include reducing the cost of the CP solve itself and optimizing when CP is applied during broad search. Such improvements would allow the framework to process more near-feasible candidates and recover feasible solutions more robustly within practical computation time.

Improving the convergence rate of the subsequent NLP is also an important future direction. In the present application, three clusters fail to converge despite containing candidates with CP-generated initial control histories, although these clusters contain relatively few representative candidates. One possible cause is that the current framework does not always provide initial guesses sufficiently close to feasible solutions of the full multi-flyby optimization problem. The broad search generates trajectory candidates on a leg-by-leg basis, and CP likewise computes the corresponding low-thrust control histories for individual legs. In contrast, the subsequent NLP optimizes the entire multi-leg trajectory simultaneously while enforcing the flyby constraints between successive legs. This difference leaves room to improve the initial guesses from the viewpoint of the complete trajectory. One possible approach is to introduce a computationally efficient convex-optimization stage between the broad search and local optimization that treats all legs simultaneously and incorporates the flyby constraints. Such an intermediate stage could screen broad-search candidates according to their likelihood of yielding feasible locally optimized solutions and generate higher-quality initial state and control histories. These improvements could increase the convergence robustness of the subsequent NLP.

\section{Conclusions} \label{sec:conclusion}

This paper proposes an integrated deterministic trajectory optimization framework for low-thrust gravity-assist mission design that automatically searches for optimal solutions over a global design space using only mission-sequence information and spacecraft specifications as inputs. The framework integrates a \textit{Star}-based broad search with Convex Pruning (CP) and local optimization using a thrust-regularized Sims-Flanagan Transcription (SFT).
This paper demonstrates the effectiveness of the proposed framework using a BepiColombo-inspired Mercury transfer mission involving nine gravity assists and low-thrust propulsion. 
The results show two main effects of this study. First, the CP in the broad-search stage discovers trajectory candidates that cannot be captured by a lightweight model alone and thereby expands the searchable design space. Second, the framework connects the initial trajectory guesses obtained from the broad search and the CP-generated initial control histories to the thrust-regularized local-optimization model, improving the convergence behavior of the subsequent trajectory optimization. Overall, the developed framework provides a new capability to optimize complex gravity-assist low-thrust trajectories in an automated manner with minimal user inputs.

\section*{Acknowledgments}
This work was conducted during R. Iijima's tenure as a visiting scholar at the School of Aeronautics and Astronautics, Purdue University. R. Iijima received support from the GP-Mech International Joint Graduate Program at Tohoku University. During the preparation of this manuscript, generative AI tools were used in part to identify more readable sentences.

\bibliography{paperpile,reference}

\end{document}